\documentclass[11pt,reqno]{amsart}

\usepackage[utf8]{inputenc}
\usepackage{amsfonts}
\usepackage{amsthm}
\usepackage{amsmath}
\usepackage{amssymb}
\usepackage{stmaryrd}
\usepackage{yfonts}
\usepackage{MnSymbol}
\usepackage{mathtools}
\usepackage{fullpage}
\usepackage{color}
\usepackage{enumerate}
\usepackage{capt-of}
\usepackage[hidelinks]{hyperref}

\usepackage{tikzscale}
\usepackage{tikz}
\usepackage{import}
\usepackage{subcaption}
\usetikzlibrary{arrows.meta,bending,positioning}
\usetikzlibrary{plotmarks}
\usepackage{adjustbox}
\usepackage{tikz-cd}
\usepackage{float}

\date{\today}

\newtheorem{teo}{Theorem}[section]

\newtheorem{coro}[teo]{Corollary}
\newtheorem{lema}[teo]{Lemma}

\newtheorem{prop}[teo]{Proposition}

\newtheorem{ex}[teo]{Example}

\theoremstyle{definition}
\newtheorem{claimInProof}{Claim}
\numberwithin{claimInProof}{teo}

\DeclareMathOperator{\w}{\omega}
\DeclareMathOperator{\C}{\mathcal{C}}
\DeclareMathOperator{\N}{\mathbb{N}}
\DeclareMathOperator{\R}{\mathbb{R}}
\DeclareMathOperator{\Homeo}{\textnormal{Homeo}}

\DeclareMathOperator{\OD}{\mathcal{DO}}

\newcommand{\cantor}{2^{\N}}

\NewDocumentCommand{\orb}{O{f} O{x}}{\mathcal{O}_{#1}(#2)}
\NewDocumentCommand{\der}{O{}}{#1^{(1)}}

\title{Countable dynamical systems with a dense orbit}
\author{Jhon Pérez\--Remolina}
\address{Escuela de Matem\'aticas, Universidad Industrial de Santander, C.P. 680001, Bucaramanga - Colombia.}
\email{prjhonf@gmail.com}

\author{Carlos Uzc\'{a}tegui-Aylwin}
\address{Escuela de Matem\'aticas, Universidad Industrial de Santander, C.P. 680001, Bucaramanga - Colombia.}
\email{cuzcatea@saber.uis.edu.co}

\date{\today}

\begin{document}

\begin{abstract}
We study countable dynamical systems with a dense orbit. Given an infinite compact metric countable space $X$, we show that $\OD(X)$, the collection of continuous functions on $X$ to itself
admitting a dense orbit, is always nonempty. We then classify the systems in $\OD(X)$
up to topological conjugacy in terms of the second Cantor--Bendixson derivative
$X^{(2)}$: if $X^{(2)} = \emptyset$ then all systems in $\OD(X)$ are mutually conjugate,
while if $X^{(2)} \neq \emptyset$ then $\OD(X)$ contains $2^{\aleph_0}$ pairwise
nonconjugate systems. The proof of the latter relies on a gluing construction for
countable dynamical systems with a dense orbit, which is of independent interest.
We also show that $\OD(X)$ is nowhere dense in $C(X,X)$ with the compact-open topology.

We also study the Ellis semigroup $E(X,f)$ of a countable dynamical system. We characterize, under the assumption that $|\w(x,f)| = 1$ for each $x \in X^{(1)}$, when the nontrivial functions of $E(X,f)$ are continuous. In particular, for the family constructed on spaces with $|X^{(2)}| \geq 2$, every nontrivial function of the Ellis semigroup is discontinuous.
\end{abstract}


\maketitle

\section{Introduction}

In this paper we study  countable dynamical systems with a dense orbit. We consider dynamical systems $(X,f)$, where $X$ is a compact metric countable space and $f:X\to X$ is continuous, and we focus on those systems for which there exists $x\in X$ such that its orbit
$$
\mathcal{O}_f(x)=\{f^n(x):n\in\mathbb{N}\}
$$
is dense in $X$. 
For a compact metric space $X$, write
\[
\OD(X) = \{f \in C(X,X) : f \text{ has a dense orbit}\}.
\]
One of the main results of this paper is the following (see Theorem \ref{OD-no-vacio}).

\medskip

\noindent\textbf{Theorem A.} If $X$ is a compact metric countable space then $\OD(X)$ is nonempty.

\medskip

Perhaps the first examples of dynamical systems on $\w^2+1$ with a dense orbit were presented in \cite{Cardinality2018}. The preceding result was motivated by a question raised in that paper.

The existence of a dense orbit is a classical manifestation of nontrivial, and often chaotic, dynamical behavior. For instance, topological transitivity is one of the defining conditions of Devaney chaos \cite{Devaney1989}. On a perfect compact metric space, it is equivalent to the existence of a dense orbit.
Another example is chaos in the sense of Martelli, which explicitly includes the existence of a dense orbit among its defining conditions \cite{Martelli}. However, these classical notions of chaos are not applicable to countable dynamical systems, since  in a countable compact metric space the set of isolated points is dense and therefore there are no transitive countable systems.  Nevertheless,  the existence of a dense orbit is a natural condition in the setting of countable spaces.

Having established the existence of dense orbits, it is natural to ask how many essentially different dynamical systems with this property can be supported by a given space. The standard notion for comparing dynamical systems in this setting is topological conjugacy. That is, two dynamical systems $(X,f)$ and $(Y,g)$ are \textit{conjugate} if there is a homeomorphism $\varphi : X \to Y$ such that $\varphi \circ f = g \circ \varphi$. Our second main result (see Theorem \ref{Teorema existencia}) shows a dichotomy in terms of the second Cantor--Bendixson derivative $X^{(2)}$ which is the set of accumulation points of $X^{(1)}$.

\medskip

\noindent\textbf{Theorem B.}  Let $X$ be a compact metric countable space. If $X^{(2)}=\varnothing$ then all systems in $\OD(X)$ are mutually conjugate. Conversely, if $X^{(2)}\neq\varnothing$, then $\OD(X)$ contains $2^{\aleph_0}$ pairwise nonconjugate systems.

\medskip

The second part of the above theorem provides an explicit construction of  $2^{\aleph_0}$ pairwise nonconjugate dynamical systems. The first step is to show this result on $\omega^2+1$. We then develop a gluing technique (see Section~\ref{section: pegado de sistemas}) for combining countable dynamical systems with dense orbits. The gluing is performed along the dense orbits of the systems being combined, and is designed so that the resulting system still has a dense orbit while retaining other properties of the systems involved. This allows us to combine the systems constructed on $\omega^2+1$ with the systems obtained in the proof of Theorem A, thereby producing the required family of pairwise nonconjugate systems.

\medskip 

Although dynamical systems on countable compact spaces have received comparatively less attention than systems on perfect spaces, they exhibit a rich variety of dynamical phenomena.

One direction that has been studied concerns the Ellis semigroup associated with a dynamical system. In particular, \cite{GarciaRodriguezTomita,  iteratesOfDN,Cardinality2018,GarciaFerreira2007} establish several properties of the Ellis semigroup arising from dynamical systems on compact metric countable spaces. More recently, Quintero and Uzcátegui \cite{QuinteroUzcategui2020} studied countable dynamical systems in the case where $f$ is a homeomorphism, obtaining results relating distality and the structure of the corresponding dynamical system.

Other dynamical properties have also been investigated in this setting. Huang and Ye  \cite{HuangAndYe} showed that there exist compact countable metric spaces admitting completely scrambled homeomorphisms. In another direction, Kato and Park \cite{KatoAndPark} characterized the compact countable metric spaces that admit expansive homeomorphisms. Bobok \cite{ChaosInCountableDS} showed that distributional chaos can occur on a countable compact metric space whenever $|X^{(2)}|>1$. A related line concerns sequence entropy on countable compact spaces. Ye and Zhang \cite{YeZhang2008} showed that homeomorphisms with positive sequence entropy exist on 
$X$ if and only if $X^{(2)}\neq \emptyset$.

A result related to Theorem A was obtained in D. Pérez's Ph.D. thesis \cite{CountableSpacesAndCountableDynamics}, where, for each ordinal $\alpha$, the author constructed a dynamical system $(X,f)$ with  Cantor-Bendixson rank  equal to $\alpha+1$ and  such that $f$ is a homeomorphism admitting a dense $\mathbb{Z}$-orbit. 

We also look at the well\--studied Ellis semigroup $E(X,f)$ associated with a dynamical system $(X,f)$ \cite{GarciaRodriguezTomita, Cardinality2018, GRC2021, QuinteroUzcategui2020}. We show that for all  systems $(X,f)$ with $|X^{(2)}| \geq 2$ constructed in this paper their Ellis semigroup is homeomorphic to $X$ and each nontrivial function of the Ellis semigroup is discontinuous.

The paper is organized as follows. Section~2 collects the necessary background on countable compact metric spaces, Cantor--Bendixson derivatives, orbits and $\omega$-limit sets. Section~3 gives explicit systems on $k\cdot\omega+1$ with a dense orbit. Section~4 constructs the families of pairwise nonconjugate systems on $\omega+1$ and $\omega^2+1$ used later. Section~5 proves  Theorem A which establishes existence for arbitrary infinite compact
metric countable spaces. Section~6 develops the gluing construction, and Section~7 combines all of the above to prove Theorem B. Finally, Section~8 studies the Ellis semigroups of the constructed systems.

\section{Preliminaries}
\label{preliminaries}

A dynamical system is a pair $(X,f)$ where $X$ is a compact metric space and $f: X \to X$ is a continuous function. In this work we use $\N = \{ 0, 1, \dots \}$. The orbit of a point $x \in X$ denoted by $\orb$ is the set $\{ f^n(x) : n \in \N \}$. The set of accumulation points of the sequence $(f^n(x))_{n \in \N}$ is called the $\w$\--limit set of $x$, which we will denote by $\w(x,f)$, more precisely,
\[
\w(x,f) = \{ y \in X : \exists \,(n_k) \subseteq \N \text{ strictly increasing such that } (f^{n_k}(x))_{k \in \N} \to y \}.
\]
A dynamical system $(X,f)$ has a dense orbit if there exists $z \in X$ such that $\orb[f][z]$ is dense in $X$. We define $\OD(X)$ as the set of continuous functions on $X$ to itself that admit a dense orbit, that is,
\[
\OD(X) = \{ f \in \C(X,X) : \exists z \in X \text{ with } \overline{\orb[f][z]} = X \}.
\]
It is well\--known that for $X$ perfect compact metric spaces, $(X,f)$ has a dense orbit exactly when it is a transitive system (that is to say, for every pair of nonempty open sets $V$ and $W$ there is $n\in \N$ such that $f^{n}[V]\cap W\neq\emptyset$).

Two dynamical systems $(X,f)$ and $(Y, g)$ are topologically conjugate if there exists a homeomorphism $\varphi: X \to Y$ such that $\varphi \circ f = g \circ \varphi$.

For any topological space $X$, we will denote by $I(X)$ the set of isolated points of $X$.
Let us denote by $X^{(1)}$ the set of all accumulation points of $X$. We call $X^{(1)}$ the {\it Cantor-Bendixson derivative} of $X$. By transfinite recursion we define the {\it$\alpha$-th Cantor-Bendixson derivative} of $X$ (denoted by $X^{(\alpha)}$) as follows: $X^{(0)} =X$, $X^{(\alpha + 1)} = (X^{(\alpha)})^{(1)}$ and 
\[
X^{(\lambda)} = \bigcap_{\alpha < \lambda} X^{(\alpha)},
\]
if $\lambda$ is a limit ordinal. The {\it Cantor-Bendixson rank} of $X$ is the first ordinal $\alpha < \omega_1$ for which $X^{(\alpha)} = X^{(\alpha + 1)}$. We will use $CB(X)$ to denote the Cantor-Bendixson rank of $X$. It is a well\--known a classical theorem that if $X$ is a compact metric countable space and $\alpha = CB(X)$, then $X^{(\alpha)} = \emptyset$. Moreover, $CB(X)$ is a successor ordinal and $I(X)$ is dense in $X$.

We denote by $\cantor$ the set of infinite binary sequences and $2^{<\omega}$ the set of finite binary sequences, that is,
\[
    2^{< \omega} = \bigcup_{n \in \N} \{ 0,1 \}^n,
\]
with the convention that $\{ 0,1 \}^0 = \{ \emptyset \}$, where $\emptyset$ denotes the empty sequence. If $s \in \{ 0,1 \}^n$, the length of $s$ is $n$ which we denoted by $|s|$. For $s = (s_0,\dots, s_n),t = (t_0, \dots, t_m) \in 2^{<\omega}$, we write  
\[
s \subseteq t \iff n \leq m \text{ and } s = (t_0, \dots, t_n).
\]
If $s = (s_0, \dots, s_{n}) \in 2^{<\omega}$ and $k \in \N$, the concatenation of $s$ and $k$ is the sequence given by $s\frown (k) = (s_0, \dots, s_n, k)$.

For any $k \in \N$ and $\alpha < \w_1$, we consider the space $\w^\alpha \cdot k + 1 = [0, \w^\alpha \cdot k]$ with the order topology.
The following is a classical well\--known result on compact metric countable spaces due to Mazurkiewicz and Sierpiński.

\begin{teo}
(Mazurkiewicz - Sierpiński)
\label{Sierpinski}
Let $X$ be a compact metric countable space. Then $X$ is homeomorphic to $\w^\alpha \cdot k + 1$ for some countable ordinal $\alpha \geq 1$ and $k \geq 1$ natural number and $CB(X)=\alpha+1$.
\end{teo}

\begin{prop}
    \label{lema: KatoPark}
    (Kato - Park \cite[Proposition 2.1]{KatoAndPark})
    Let $X$ and $Y$ be two compact metric countable spaces. If $CB(X) = CB(Y) = \alpha + 1$ and $X^{(\alpha)}$ is homeomorphic to $Y^{(\alpha)}$, then $X$ is homeomorphic to $Y$.
\end{prop}

\begin{coro}
    \label{Corolario: homeomorfismo de puntos limites es suficiente}
    Let $X$ and $Y$ be two compact metric countable spaces. If $X^{(1)}$ is homeomorphic to $Y^{(1)}$, then $X$ and $Y$ are homeomorphic.
\end{coro}

Whenever we consider the union of two compact metric countable spaces, we regard them as subspaces of $[0,2]\times[0,1]$. Thus, for compact metric countable spaces $X_1$ and $X_2$, the union $X_1\cup X_2$ is endowed with the subspace topology inherited from $[0,2]\times[0,1]$.

\begin{prop}
    \label{union homeomorfa a uniendo}
    Let $X_1, X_2$ be two compact metric countable spaces and let $\lambda$ be such that $ CB(X_2)=\lambda + 1$. If $CB(X_1) \geq CB(X_2)$ and $X_2^{(\lambda)} \subseteq X_1^{(\lambda)}$, then $X_1 \cup X_2$ and $X_1$ are homeomorphic. 
\end{prop}
\begin{proof}
    It is clear that $X: = X_1 \cup X_2$ is a compact metric countable space. Note that $X^{(\alpha)} = X_1^{(\alpha)} \cup X_2^{(\alpha)}$ for all $\alpha < \w_1$. Let $\alpha + 1 = CB(X_1)$. Since $X_2^{(\lambda)} \subseteq X_1^{(\lambda)}$ and $\lambda \leq \alpha$ we get $X_2^{(\alpha)} \subseteq X_1^{(\alpha)}$. Therefore $X^{(\alpha)} = X_1^{(\alpha)}$. Thus, the conclusion follows from Proposition \ref{lema: KatoPark}.
\end{proof}

\subsection{Topological dynamics}

\begin{prop}
    \label{Proposicion: funcion con orbita densa}
    Let $(X,f)$ be a countable dynamical system and $z \in X$. The orbit of $z$ is dense if and only if $\orb[f][z] = I(X)$. 
\end{prop}

\begin{proof} It is well\--known that in a countable topological space $I(X)$ is dense  in $X$ and moreover any dense set $D$ must contain $I(X)$. In particular,  whenever an orbit is dense, it must be infinite. On the other hand, $\orb[f][f^m(z)]\cap I(X)\neq \emptyset$ for all $m$. By the continuity of $f$ and using that the orbit of $z$ is dense, if for some $t$ and $m$, $f^m(t)\in I(X)$, then $t\in I(X)$.
\end{proof}

The following are well\--known facts about $\w(x,f)$. 

\begin{prop}
\label{Lema: resutlado clasico w(x,f) contenido en U}
Let $(X,f)$  be a dynamical system and  $x \in X$. Then 
\begin{enumerate}[(i)]
\item $\w(x,f)$ is a nonempty closed set.
\item $f(\w(x,f)) = \w(x,f)$. 
\item  If $\w(x,f)$ is finite, then $\w(x,f)$ is an orbit.
\item 
\label{prop: iteradas en un abierto}
If $U \subseteq X$ is an open set and  $\w(x,f) \subseteq U$, then there is $K \in \N$ such that $f^n(x) \in U$ for each $n \geq K$.
\end{enumerate}    
\end{prop}

In the following lemma we collect some facts about the dynamics when $\w(x,f)$ is finite. Perhaps they are known but we could not find a proper reference for them. 

For $k, r \in \N$ we denote
\[
k\N + r = \{ kn + r : n \in \N \}.
\]

As noted above, when $\w(x,f)$ is finite, it is necessarily a periodic orbit. Thus we can enumerate $\w(x,f)=\{x_0,\ldots, x_{k-1}\}$ such that 

\begin{enumerate}
\item[(a)]  $f(x_i)=x_{i+1}$ for $0\leq i<k-1$ and $f(x_{k-1})=x_0$. 
\end{enumerate}
Moreover  we can easily find open sets $V_i$, for $0\leq i<k$,  such that 

\begin{enumerate}
\item[(b)] $\overline{V_i} \cap \overline{V_j} = \emptyset$ for $i \neq j$.

\item[(c)] $x_i\in V_i$ for all $0\leq i<k$. 

\end{enumerate}

\begin{lema}
\label{lema de sincronizacion}
Let $(X,f)$ be a dynamical system, $x \in X$ with  $k = |\w(x,f)|$ and let $V_0, \dots, V_{k - 1} \subseteq X$ be open sets satisfying (a), (b) and (c).  Then the following holds:

\begin{enumerate}[(i)]
\item There exists $L \in \N$ such that for each $n \geq L$ then $f^n(x) \in V_i$ for some $i < k$.

\item There is $L$ such that for each $n \geq L$, if $f^n(x) \in V_i$ for some $i$ then $f^{n+ 1}(x) \in V_{i + 1 (\textnormal{mod } k)}$.
        
\item  
\label{lema-item: omega l.s. particiona}
For each $0 \leq i < k$ there is $0 \leq j < k$ such that $\{ n \in \N: f^n(x) \in V_i \} \setminus (k\N + j)$ is finite.

\item There is $m \in \N$ such that
\[
d(f^{m + j}(x), f^{j}(x_0)) \leq \textnormal{diam}(V_{j (\textnormal{mod } k)}) \ \text{ for all } j \in \N.
\]
\end{enumerate}
\end{lema}

\begin{proof}
\textit{(i)}. By Proposition \ref{Lema: resutlado clasico w(x,f) contenido en U}\eqref{prop: iteradas en un abierto}  and (b), there is $L \in \N$ such that for all $n \geq L$, there is a unique $0 \leq i(n) < k$ such that $f^n(x) \in V_{i(n)}$.  

(ii) Note that this is equivalent to the following: There exists $L^{\prime}\geq L$ such that for each $n \geq L^{\prime}$,
\[
f^n(x) \in V_{i(n)} \  \Longrightarrow \ f^{n+1}(x) \in V_{i(n) + 1 \ (\text{mod } k)}.
\]
In fact, towards a contradiction, suppose there is a strictly increasing  sequence $(n_j)_{j \in \N}$ with each $n_j\geq L$ such that
\[
            f^{n_j}(x) \in V_{i(n_j)} \; \text{ and } \; f^{n_j + 1}(x) \notin V_{i(n_j) + 1 \ (\text{mod } k)}.
\]
We can assume w.l.o.g. there is $l<k$ such that $n_j \equiv l \ (\text{mod } k)$ for all $j$ and $(f^{n_j}(x))_{j \in \N}$ is a convergent sequence (by the compactness of $X$). On the other hand, since $\w(x,f) \cap V_l = \{ x_l \}$ (recall the enumeration of $\w(x,f)$ and condition (a)) we have $(f^{n_j}(x))_{j \in \N} \to f^l(x_l)$. By continuity, $(f^{n_j + 1}(x))_{j \in \N}$ converges to $f^{l + 1}(x_l)$ and by (a), $f^{l + 1}(x_l) \in V_{l + 1 \ (\text{mod } k)}$. Thus $f^{n_j + 1}(x) \in V_{l + 1\ (\text{mod } k)}$ for all sufficiently large $j \in \N$, which is a contradiction.
       
(iii) It follows from (ii). 

(iv) By (iii), for each $i < k$ there is $j_i$ such that $\{ n \in \N : f^n(x) \in V_i \} \setminus (k\N + j_i)$ is finite. Hence there is  $L''$ such that 
$\{ n \in \N:n>L''\;\& \; f^n(x) \in V_i \} \subseteq k\N + j_i$, for all $0\leq i<k$.
Take $m=\min \{ n \in \N:n>L''\;\& \; f^n(x) \in V_0 \}$.
\end{proof}

Although part (iv) of the preceding lemma will not be used, it illustrates  a sort of synchronization occurring on the orbit. 

\medskip 
    
\begin{lema}
\label{caracterizacion convergencia en conjuntos omega limites}
Let $(X,f)$ be a dynamical system, $x \in X$,  $|\w(x,f)| =k$ and  $A \subseteq \N$ be an infinite subset. Then 
\[
(f^n(x))_{n \in A} \text{ converges } \iff \exists j<k \ ( A \setminus (k\N + j) \text{ is finite} ).
\]
\end{lema}

\begin{proof}
Let $\w(x,f)=\{x_0, \ldots, x_{k-1}\}$ and $V_i$, for $0\leq i< k$ be as in Lemma \ref{lema de sincronizacion}. 
Let us suppose that $(f^n(x))_{n \in A}$ converges to $x_i$.  The conclusion follows directly from Lemma \ref{lema de sincronizacion}\eqref{lema-item: omega l.s. particiona}.

Conversely, if $A$ is almost contained in $k\N + j$, then by Lemma \ref{lema de sincronizacion}\eqref{lema-item: omega l.s. particiona}, necessarily $f^n(x) \in V_i$ (for some $i < k$) for almost all $n \in A$. Thus 
$(f^n(x))_{n \in A}$ converges to $x_i$. 
\end{proof}

The following lemmas are straightforward.
\begin{lema}
    \label{Lema: preimagen conjugacion}
    Let $(X,f)$ and $(Y, g)$ be two dynamical systems. If $\varphi: X \to Y$ is a conjugacy map between $(X,f)$ and $(Y,g)$, then 
    \[
    \varphi[f^{-n}[x]] = g^{-n}[\varphi(x)],
    \]
    for all $n \in \N$ and for each $x \in X$.
\end{lema}

\begin{lema}
   \label{OrbitaDensaBajoHomeomorfismos}
   Let $X$ and $Y$ be two compact metric spaces, and let $\varphi: Y \to X$ be a homeomorphism. If $(Y,g)$ has a dense orbit, then so does $(X, \varphi \circ g \circ \varphi^{-1})$.
\end{lema}

\section{The case \texorpdfstring{$\w \cdot k + 1$}{wk + 1}}
\label{Section: case kw + 1}

\begin{figure}[H]
\centering
\begin{tikzpicture}[
    every label/.append style={font=\footnotesize},
    scale=1.24
  ]

    \tikzstyle{point}=[draw,shape=circle, inner sep=.75pt, fill]
    \tikzstyle{arrow} = [-{Classical TikZ Rightarrow[scale=1.1, sep=1pt, length=1mm]}]
  
    \newcommand{\drawBlock}[2]{
        \def\listDistance{0,1.5,2.7,3.7,4.6,5.3}
        
        \pgfmathtruncatemacro{\nmax}{#2}
        \pgfmathsetmacro{\xlim}{#1}
        \pgfmathsetmacro{\xaux}{#1}
        \pgfmathsetmacro{\dlast}{0}
        
        \foreach \d [count=\m from 0] in \listDistance {
            \ifnum\m > \nmax
                \relax
            \else
                \pgfmathsetmacro{\x}{\xlim + \d}
                \node[point, label={[label distance=.5mm]270:${\m}$}] (\m) at (\x, 0) {};
                \pgfmathparse{\x}
                \xdef\xaux{\pgfmathresult}
                \pgfmathsetmacro{\dlast}{\d}
            \fi
        }
        \pgfmathsetmacro\dotsp{\xaux + 0.3}
        \node[] (suspensive) at (\dotsp, 0) {$\cdots$};

        \pgfmathsetmacro\limp{\xaux + 0.6}
        \node[point, label={[label distance=.5mm]270:$\w$}, inner sep=1.1pt] (limit) at (\limp, 0) {};
    }

    \drawBlock{-2}{5}

    \foreach \from/\to  in {
    0/1,
    1/2,
    2/3,
    3/4,
    4/5}
    {
        \draw[arrow] (\from) to [bend left] (\to);
    }

    \draw[arrow] (limit) to [loop] (limit);
\end{tikzpicture}
\caption{Dynamical system $(\w + 1,f)$ with dense orbit}
\label{fig:dn-on-w}
\end{figure}
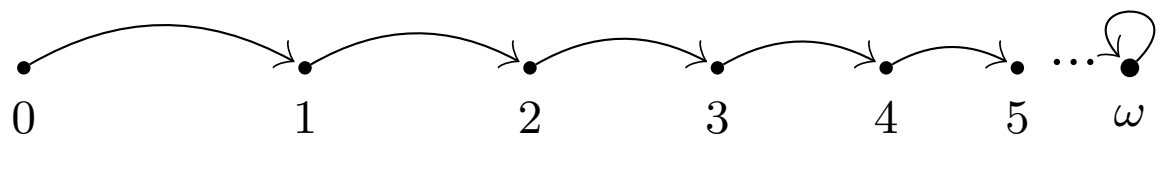

In this section we present some elementary examples of dynamical systems on $\w \cdot k + 1$. Such examples are straightforward to construct when $k = 1$. Recall that $\w + 1 = \N \cup \{ \w \}$. The space $\omega+1$  has exactly one limit point and any function $f: \N \to \N$ such that $\N$ is the orbit of a point provides an example of a dynamical system on $\omega+1$ with a dense orbit by declaring $f(\omega)=\omega$. Figure \ref{fig:dn-on-w} illustrates one such example.

For $k > 1$, the construction remains fairly simple. As an illustration, we present a system on $\w\cdot 3 + 1$.

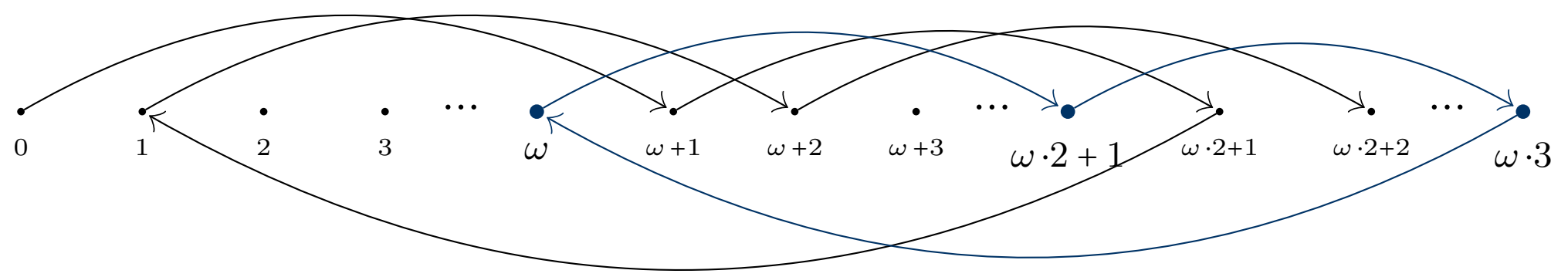
\begin{figure}[H]
\centering
\begin{tikzpicture}[
    every label/.append style={font=\footnotesize},
    scale=1.3
  ]

    \tikzstyle{point}=[draw,shape=circle, inner sep=.5pt, fill]
    \tikzstyle{arrow} = [-{Classical TikZ Rightarrow[scale=1.1, sep=1pt, length=1mm]}]

    \definecolor{azulejo}{RGB}{0,51,102}

    \foreach \p [count=\m from 0] in {0, 0.8, 1.6, 2.4}{
        \node[point, label={[label distance=0.7mm]270:$ \scriptstyle \m$}] (1\m) at (\p, 0) {};
    }

    \node[] (suspensive1) at (2.9, 0) {$\cdots$};

    \node[point, color=azulejo, label={[label distance=0.7mm]270:\footnotesize{$\w$}}, inner sep=1.1pt] (limit1) at (3.4, 0) {};

    \foreach \p [count=\m from 1] in {4.3, 5.1, 5.9}{
        \node[point, label={[label distance=0.7mm]270:$\scriptstyle \w +  \m$}] (2\m) at (\p, 0) {};
    }

    \node[] (suspensive2) at (6.4, 0) {$\cdots$};

    \node[point, color=azulejo, label={[label distance=0.7mm]270:\footnotesize{$\w\cdot 2 + 1$}}, inner sep=1.1pt] (limit2) at (6.9, 0) {};

    \foreach \p [count=\m from 1] in {7.9, 8.9}{
        \node[point, label={[label distance=0.7mm]270:$\scriptstyle \w\cdot 2 +  \m$}] (3\m) at (\p, 0) {};
    }

    \node[] (suspensive2) at (9.4, 0) {$\cdots$};
    
    \node[point, color=azulejo, label={[label distance=0.7mm]270:\footnotesize {$\w\cdot3$}}, inner sep=1.1pt] (limit3) at (9.9, 0) {};

    \foreach \from/\to  in {
    10/21,
    21/31,
    31/11,
    11/22,
    22/32}
    {
        \draw[arrow] (\from) to [bend left] (\to);
    }

    \foreach \from/\to  in {
    limit1/limit2,
    limit2/limit3,
    limit3/limit1}
    {
        \draw[draw=azulejo,arrow] (\from) to [bend left] (\to);
    }
  
\end{tikzpicture}
\caption{Dynamical system with dense orbit on $\w\cdot 3 + 1$}
\end{figure}

Consider the sets $M_i = \{ \w \cdot i + \ m: m \in \N \}$, for $i = 0,1,2$. The orbit of the point $0$ can then be described as one that visits each set $M_i$ successively, in a fixed cyclic order. 

The previous construction extends naturally to yield dynamical systems with a dense orbit on $k \cdot \w + 1$ for arbitrary $k > 1$. More precisely, for $k > 1$ consider the function $f :\w \cdot k + 1 \to \w \cdot k + 1$ given by
\[
    f(x) = \begin{cases}
            \w + (m + 1) & \text{ if } x = m,\\
            \w \cdot (n + 1) + \ m  & \text{ if } x = \w\cdot n + \ m, \text{ } 0 < n < k -1 \text{ and } m \geq 1,\\
            m + 1 & \text{ if } x = \w \cdot (k-1) + \ m \text{ and } m \geq 1,\\
            \w \cdot (n + 1) & \text{ if } x = \w \cdot n \text{ and } 0 < n < k,\\
            \w & \text{ if } x = \w\cdot k.
        \end{cases}
\]
The continuity of $f$ can be verified using sequences. For instance, let $0<n<k-1$, then $\omega \cdot n + i \to \omega\cdot (n + 1) + i$ and,
$$
f(\omega\cdot n +i)=\omega\cdot (n + 1) +i \longrightarrow \omega \cdot (n + 2) =f(\omega\cdot(n + 1)).
$$
Also observe that the orbit of $0$ is $I(X)$. Consequently $(X,f)$ has a dense orbit (see Proposition \ref{Proposicion: funcion con orbita densa}).

While the construction of dynamical systems on spaces of the form $\w\cdot k + 1$ presents no  difficulty, the situation becomes more delicate once the Cantor--Bendixson rank of the space increases. Nevertheless, in \cite{Cardinality2018} the authors construct examples of such systems on $\w^2 + 1$ and  $\w^3 + 1$.

\section{Construction of dynamical systems  on \texorpdfstring{$\w^2 + 1$}{w2 + 1} with a dense orbit}

In this section  we show how to construct a family of continuum many pairwise nonconjugate dynamical systems with a dense orbit on $\omega^2+1$.  This establishes one case of the main theorem. The construction on $\omega^2+1$ also serves as a guide for the general case of a compact metric  countable space. Two examples of  dynamical systems on $\w^2+1$ with a dense orbit were presented in \cite[Examples 4.2 and 4.3]{Cardinality2018}.

\subsection{Some dynamical systems on \texorpdfstring{$\w + 1$}{w + 1}}

We present the first step towards the main theorem of this section, a family of continuum many pairwise nonconjugate dynamical systems (without dense orbits) on a convergent sequence.

Given $\alpha \in \cantor$, we define a sequence $(l_n^\alpha)_{n \in \N}$ of natural numbers by $l_0^\alpha = 1$ and $l_{n + 1}^\alpha = (1 + \alpha(n))\cdot l_n^\alpha$.
Now, we will define a tree $T_\alpha$ as follows: 
\[
    T_\alpha = \{ s \in 2^{<\w} :  s(i) \leq \alpha(i) \text{ for each } i < |s|\}.
\]
Intuitively $T_\alpha$ is a tree in which each node at the $n$-th level has one or two immediate successors, depending on whether $\alpha(n)$ is $0$ or $1$. We define the levels of the tree $T_\alpha$ as $\widehat{L}_n^\alpha = \{ s \in T_\alpha : |s| = n \}$.

\begin{lema}
    Let $\alpha \in \cantor$, and let $T_\alpha$ be the tree defined above. Then $|\widehat{L}_n^\alpha| = l_n^\alpha$ for $n \in \N$. 
\end{lema}
\begin{proof}
    We proceed by induction on $n$. Clearly $|\widehat{L}_0^\alpha| = l_0$. Assume that $|\widehat{L}_k^\alpha| = l_k^\alpha$. Let us observe that $\widehat{L}_{k+1}^\alpha =  \{ s \frown (i) : s \in \widehat{L}_k^{\alpha} \, \text{ and } \, i \leq \alpha(k) \}$. Hence $|\widehat{L}_{k+1}^\alpha| = (1 + \alpha(k)) \cdot |\widehat{L}_k^\alpha|$, thus, by induction hypothesis we conclude $|\widehat{L}_{k + 1}^\alpha| = (1 + \alpha(k)) \cdot l_k^\alpha = l_{k + 1}^\alpha$.
\end{proof}

Let $\{ L_n^\alpha \}_{n \in \N}$ be a partition of $\N$ such that $|L_n^\alpha| = l_n^\alpha$ and $L_0^\alpha = \{ 0 \}$.
By the previous lemma, for each $\alpha \in \cantor$ fix a bijection $j_\alpha : \bigcup_{n \in \N} L_n^\alpha \to \bigcup_{n \in \N} \widehat{L}_n^\alpha$ such that $j_\alpha[L_n^\alpha] = \widehat{L}_{n}^\alpha$ for every $n \in \N$.

\begin{figure}[H]
\centering
\begin{tikzpicture}[
    every label/.append style={font=\footnotesize},
    scale=0.8
  ]

    \definecolor{verdePersonalizado}{RGB}{91,226,1}
    \definecolor{azulPersonalizado}{RGB}{82,100,255}
     \tikzstyle{arrow} = [-{Classical TikZ Rightarrow[scale=1.1, sep=1pt, length=1mm]}]
    \tikzstyle{arrowRed} = [-{Classical TikZ Rightarrow[scale=1.1, sep=1.8pt, length=2mm, color=red]}]

    \tikzstyle{point}=[draw,shape=circle, inner sep=.75pt, fill]

    \node[] (l0) at (-5, -.4) {$\scriptstyle 0$};

    \node[] (l4) at (-7, -1.2) {$\scriptstyle 1$};

    \node[] (l6) at (-3, -1.2) {$\scriptstyle 2$};

    \node[] (l8) at (-7.8, -2.5) {$\scriptstyle 3$};
    \node[] (l10) at (-6.2, -2.5) {$ \scriptstyle 4$};

    \node[] (l12) at (-3.8, -2.5) {$\scriptstyle 5$};
    \node[] (l14) at (-2.2, -2.5) {$\scriptstyle 6$};

    \node[] (l16) at (-7.8, -3.8) {$\scriptstyle 7$};
    \node[] (l18) at (-6.2, -3.8) {$\scriptstyle 8$};

    \node[] (l20) at (-3.8, -3.8) {$\scriptstyle 9$};
    \node[] (l22) at (-2.2, -3.8) {$\scriptstyle 10$};

    \node[point] (n1) at (-8.2, -4.7) {};
    \node[point] (n2) at (-7.4, -4.7) {};

    \node[point] (n3) at (-6.6, -4.7) {};
    \node[point] (n4) at (-5.8, -4.7) {};

    \node[point] (n5) at (-4.2, -4.7) {};
    \node[point] (n6) at (-3.4, -4.7) {};
    
    \node[point] (n7) at (-2.6, -4.7) {};
    \node[point] (n8) at (-1.8, -4.7) {};
    
    \node[point] (n9) at (-8.2, -5.3) {};
    \node[point] (n10) at (-7.4, -5.3) {};
    \node[point] (n11) at (-6.6, -5.3) {};
    \node[point] (n12) at (-5.8, -5.3) {};
    \node[point] (n13) at (-4.2, -5.3) {};
    \node[point] (n14) at (-3.4, -5.3) {};
    \node[point] (n15) at (-2.6, -5.3) {};
    \node[point] (n16) at (-1.8, -5.3) {};

    \foreach \from/\to in {
        l4/l0,
        l6/l0,
        l8/l4,
        l10/l4,
        l12/l6,
        l14/l6,
        l16/l8,
        l18/l10,
        l20/l12,
        l22/l14,
        n1/l16,
        n2/l16,
        n3/l18,
        n4/l18,
        n5/l20,
        n6/l20,
        n7/l22,
        n8/l22,
        n9/n1,
        n10/n2,
        n11/n3,
        n12/n4,
        n13/n5,
        n14/n6,
        n15/n7,
        n16/n8}
    {
        \draw[arrow] (\from) to (\to);
    }

    \foreach \position in {-7.8,-6.2,-3.8,-2.2}
        \node (\position) at (\position, -5.6) {\rotatebox{90}{$\cdots$}};

\end{tikzpicture}
\caption{Hasse diagram of $\preceq_\alpha$ for $\alpha=(1,1,0,1,0,0,\dots)$.}
    \label{Figure: Hasse diagram}
\end{figure}
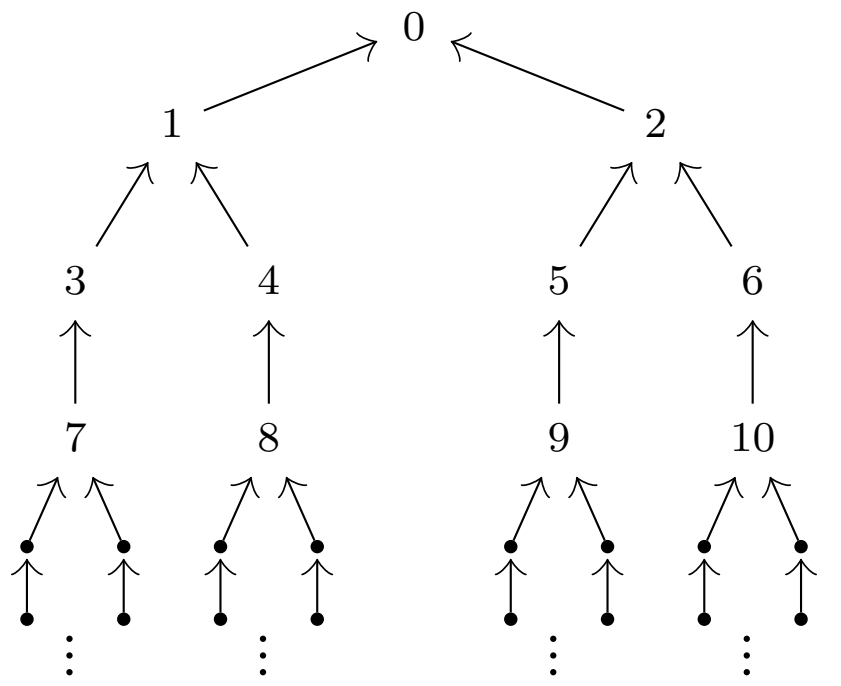

The tree structure of $T_\alpha$ induces a partial order on $\N$ by the reverse extension order. More precisely, define $\preceq_\alpha$ on $\N$ by
\[
n \preceq_\alpha m
\quad\Longleftrightarrow\quad
j_\alpha(m)\subseteq j_\alpha(n),
\]
where $\subseteq$ denotes the extension order on finite sequences. Figure \ref{Figure: Hasse diagram} is the Hasse diagram of $\preceq_\alpha$ where $\alpha = (1,1,0,1,0,0,\dots)$ and $L_0 = \{ 0 \}, L_1 = \{ 1,2 \}$, $L_2 = \{ 3,4,5,6 \}$, etc.

\noindent

We now define a function $s_\alpha : \w + 1 \to \w + 1$ associated with $\preceq_\alpha$. Let
\begin{equation}
    \label{eq: definicion s_alfa}
    s_\alpha(x) = \begin{cases}
        \w & \text{ if } x = \w \text{ or } x = 0,\\
        j_\alpha^{-1}(s)  & \text{ if } x = j_\alpha^{-1}(s \frown (k)) \text{ for some } k \in \{  0,1 \}.
    \end{cases}
\end{equation}

$s_\alpha$ maps each node of the tree to its immediate predecessor, while both the root $0$ and the point $\omega$ are mapped to $\omega$. Thus, in some sense, $s_\alpha$ can be regarded as the successor function associated with the order $\preceq_\alpha$. For instance, according to the Hasse diagram in Figure \ref{Figure: Hasse diagram}, we have $s_\alpha(7)=3$ and $s_\alpha(5)=2$.

The following lemma summarizes the basic properties of $s_\alpha$.

\begin{lema}
    \label{Lema: props de s_alfa}
    Let $\alpha \in \cantor$ and let $s_\alpha : \w + 1 \to \w + 1$ be defined as in \eqref{eq: definicion s_alfa}. Then:
    \begin{enumerate}[(i)]
    
        \item 
        \label{Lema: props de s_alfa - 1}
        $s_\alpha^{-1}[L_i^\alpha] = L_{i+1}^\alpha$ for every $i \in \N$.
        
        \item 
        \label{Lema: props de s_alfa - 2}
        For every $k \in \N$, 
        $$
        s_\alpha^{-(k + 2)}[\w] = \{ \w \} \cup \bigcup_{i = 0}^{k + 1} L_i^\alpha.
        $$
        
        \item 
        \label{Lema: props de s_alfa - 3}
        If $\alpha \neq \beta$ and $k = \min \{ n \in \N : \alpha(n) \neq \beta(n) \}$, then
        \[
        \big\vert s_\alpha^{-(k + 2)}[\w]\big | \neq  \big| s_\beta^{-(k + 2)}[\w] \big|.
        \]
        
        \item 
        \label{Lema: props de s_alfa - 4}
        $(\w + 1, s_\alpha)$ is a dynamical system for each $\alpha \in \cantor$.

        \item 
        \label{Lema: props de s_alfa - 5}
        $s_\alpha$ and $s_\beta$ are not topologically conjugate for all $\alpha \neq \beta$. 
    \end{enumerate}
\end{lema}

\begin{proof}
    \textit{(i)} Let $y \in L_{i +1}^\alpha$. Let us observe that,
    \[
        y \in L_{i + 1}^\alpha \iff \text{there is } s \in \widehat{L}_{i}^\alpha \text{ such that } y = j_\alpha^{-1}(s \frown (k)) \text{ for some } k \in \{ 0,1 \}.
    \]

    From the above, $y \in L_{i + 1}^\alpha$ if and only if $s_\alpha(y) = j_\alpha^{-1}(s)\in L_{i}^\alpha$.

    \textit{(ii)} We argue by induction on $k$. For $k = 0$, since $L_0^\alpha = \{ 0 \}$ and $s_\alpha^{-1}[\w] = \{ \w \} \cup L_0^\alpha$, then $s_\alpha^{-2}[\w] = \{ \w \} \cup L_0^\alpha \cup L_1^\alpha$ by \textit{(i)}. Assume that equality holds for $k$. Then, by the induction hypothesis
    \[
    s_\alpha^{-(k+3)}[\w] = s_\alpha^{-1}
    [ \w ] \cup \bigcup_{i = 0}^{k+1} s_\alpha^{-1}[L_{i}^\alpha].
    \]
    Using part \eqref{Lema: props de s_alfa - 1} we conclude 
    \[
        s_\alpha^{-(k+3)}[\w] = \{ \w, 0 \} \cup \bigcup_{i = 1}^{k + 2} L_i^\alpha.
    \]
    Hence, since $\{ 0 \} = L_0^\alpha$ we have completed the induction step.

    \textit{(iii)} By \eqref{Lema: props de s_alfa - 2} we get 
    \[
        \big\vert s_\alpha^{-(k+2)}[\w] \big\vert = 1 + \sum_{i = 0}^{k+1} l_i^\alpha,
    \]
    for each $\alpha \in \cantor$ and $k \in \N$. Thus, if $k = \min \{ n \in \N : \alpha(n) \neq \beta(n) \}$ then $l_i^\alpha = l_i^\beta$ for $i \leq k$ and $l_{k + 1}^\alpha \neq l_{k + 1}^\beta$ by the definition of the numbers $l_i^\alpha$. This proves \eqref{Lema: props de s_alfa - 3}.
    
    \textit{(iv)} We show that $s_\alpha$ is continuous. Since every point of $\N$ is isolated, it suffices to verify continuity at $\w$. Let $U$ be a neighborhood of $\w$. Then $(\w + 1) \setminus U$ is finite. Since $s_\alpha$ is finite-to-one, the set $s_\alpha^{-1}[(\w + 1) \setminus U]$ is finite. As $s_\alpha(\w) = \w$ we obtain $s_\alpha^{-1}[U]$ is a neighborhood of $\w$, proving that $s_\alpha$ is continuous.

    \textit{(v)} Towards a contradiction we assume that $s_\alpha$ and $s_\beta$ are topologically conjugate. Let $\varphi \in \Homeo(\w + 1)$ be such that $\varphi \circ s_\alpha = s_\beta \circ \varphi$. Let us observe that we must have $\varphi(\w) = \w$ since $\w$ is the only limit point of $\w + 1$.
    On the other hand, let $k = \min \{ i \in \N: \alpha(i) \neq \beta(i) \}$. 
    Therefore, using that  $\varphi(\w) = \w$, $\varphi$ is a bijection and \eqref{Lema: props de s_alfa - 3} we get 
    \[
    \big| \varphi[s_\alpha^{-(k + 2)}[\w]] \big| 
    \neq  
    \big| s_\beta^{-(k + 2)}[\varphi(\w)] \big|.
    \]
    This contradicts Lemma \ref{Lema: preimagen conjugacion}.
\end{proof}

\subsection{Some dynamical systems on \texorpdfstring{$\w^2 + 1$}{w2 + 1}}

We will construct dynamical systems on $\w^2 + 1$ in two steps. First we are going to copy the dynamical systems  $(\w + 1, s_\alpha)$ constructed in the previous section, in a system with  phase space a subset of the real interval $[1,2]$. Second, we shall extend them systems to  a dynamical systems with a dense orbit on spaces homeomorphic to $\w^2 + 1$. 

For each $n \in \N$ define $w_n = 1 + 1/(n+1)$ and let
\[
    W = \{ w_n : n  \in \N \} \cup \{ 1 \}.
\]
$W$ (with the usual topology) is homeomorphic to $\w + 1$ and a homeomorphism from $\w + 1$ to $W$ is given by $n \mapsto w_n$ and  $\w \mapsto 1$. Thus, we consider the function $s_\alpha$ defined on $W$ instead of $\w + 1$ and  we will regard $\preceq_\alpha$ as defined on $W$ (hence $1$ is the maximum on the order $\preceq_\alpha$).

For every $x \in W$, set 
\begin{equation}
\label{segmento-alpha}
[x]_{\alpha} = \{ w \in W : x \preceq_\alpha w \}.    
\end{equation}

\begin{lema}
\label{propiedades de G_alpha}
The following assertions hold:
\begin{enumerate}[(i)]
\item
\label{imagen de G_alpha}
$s_\alpha[[w_i]_{\alpha} \setminus \{1,w_0\}] = [w_i]_{\alpha} \setminus \{ 1, w_i \}$ for each $i \in \N$ and $\alpha \in \cantor$.

\item 
\label{conjugacion G_alpha 2}
The dynamical systems defined on $W \times \{0\}$ by $(x,0) \mapsto (s_\alpha(x),0)$  are pairwise nonconjugate.
\end{enumerate}
\end{lema}
\begin{proof}
    \textit{(i)} This follows immediately from the fact that $s_\alpha$ is the successor function associated with $\preceq_\alpha$ on $W \setminus \{1,w_0 \}$.

    \textit{(ii)} The conclusion is an immediate consequence of Lemma \ref{Lema: props de s_alfa}\eqref{Lema: props de s_alfa - 5}.
\end{proof}

We continue by proving the main result of this section.

\begin{lema}
    \label{lema: Existencia sistemas dinamicos en w2 + 1}
    For each $\alpha \in \cantor$ there exists a countable dynamical system $(Y_\alpha, f_\alpha)$ such that
    \begin{enumerate}[(i)]
        \item 
        \label{lema-item: Yalpha homeo a w2 + 1}
        $Y_\alpha$ is homeomorphic to $\w^2 + 1$.
        
        \item
        \label{lema-item: condicion sobre Yalpha}
        $Y_\alpha \subseteq [1,2]\times[0,1]$, $(Y_\alpha)^{(1)} \subseteq [1,2]\times\{ 0 \}$ and $(Y_\alpha)^{(2)} = \{ (1,0)\}$.
        
        \item
        \label{lema-item: condicion M es cerrado fi}
        $M = W \times \{ 0 \}$ is closed and strongly invariant under $f_\alpha$.
        
        \item 
        \label{lema-item: condicion punto fijo falpha}
        $(1,0)$ is the only  fixed point of $(Y_\alpha, f_\alpha)$. Furthermore,
        \[
            \bigcup_{n \in \N} f_\alpha^{-n}[(1,0)] = M.
        \]
        
        \item 
        \label{lema-item: condicion punto con orbita densa}
        The point $(1,1)$ has dense orbit under $f_\alpha$.
        \item 
        \label{lema-item: condicion falpha y fbeta no conjugados}
        $f_\alpha$ and $f_\beta$ are not conjugate for $\alpha \neq \beta$. Moreover the restrictions $f_\alpha \restriction M$ and $f_\beta \restriction M$ are also not conjugate.
    \end{enumerate}
\end{lema}
\begin{proof}
    Recall $W  =\{ 1 + 1/(n + 1) : n \in \N \} \cup \{  1\} =\{w_n: n \in \N\} \cup \{ 1 \}$. For each $\alpha \in \cantor$, let $M^\alpha_{2i + 1} = \{ 1 \}$ and $M^\alpha_{2i} = [w_i]_{\alpha}$ for all $i \in \N$ (recall that $[x]_{\alpha}$ was defined in \eqref{segmento-alpha}).  Set 
    \[
    Y_\alpha = W \times \{ 0 \} \cup \bigcup_{i \in \N}  M_i^\alpha \times \{ \phi(i) \},
    \]
    where $\phi: \N \to [0,1]$ is given by $\phi(i) = \frac{1}{i + 1}$. To simplify notation, we will omit the superscript $\alpha$ and simply write $M_i$ instead of $M_i^\alpha$.
    
    With the usual topology inherited from $[1,2] \times [0,1]$, $Y_\alpha$  is a countable compact metric space. 
    Observe that $M_i$ is finite for each $i \in \N$ and the set of isolated points of $Y_\alpha$ is 
    \[ \bigcup_{i \in \N} M_i \times \{ \phi(i)\}. 
    \]
    This implies that $(Y_\alpha)^{(1)} = W \times \{ 0 \} \subseteq [1,2] \times [0,1]$. Furthermore, $(Y_\alpha)^{(2)} = \{ (1,0) \}$. The above proves \eqref{lema-item: condicion sobre Yalpha} and \eqref{lema-item: Yalpha homeo a w2 + 1}.

    Now, consider $f_\alpha: Y_\alpha \to Y_\alpha$ given by $f_\alpha(x,0) = (s_\alpha(x),0)$ for $(x, 0) \in W \times \{0\}$ and if $(x, \phi(i)) \in M_i \times \{ \phi(i) \}$ for some $i \in \N$, let 
    \[
        f_\alpha(x, \phi(i)) = \begin{cases}
            (1, \phi(i + 1)) & \text{ if } i \text{ is odd},\\[4pt]
            (s_\alpha(x), \phi(i)) & \text{ if } i \text{ is even and } x \neq 1,2,\\[4pt]
            (w_{i/2}, \phi(i)) & \text{ if } i \text{ is even and } x=1,\\[4pt]
            (1, \phi(i + 1)) & \text{ if } i \text{ is even and } x = 2.
        \end{cases}
    \]
    Figure \ref{fig: sistema dinamico f_alfa} we illustrate $f_\alpha$ for $\alpha=(0,0,0,\ldots)$.

    \begin{figure}[H]
\centering
\begin{tikzpicture}[
    every label/.append style={font=\footnotesize},
    scale=1.1
  ]

    \tikzstyle{point}=[draw,shape=circle, inner sep=.9pt, fill]
    
    \tikzstyle{ipoint}=[draw,shape=circle, inner sep=.5pt, fill]

    \definecolor{azulejo}{RGB}{0,51,102}
  
    \tikzstyle{arrow} = [
    -{Classical TikZ Rightarrow[
        scale=1.4,
        sep=2pt,
        length=1mm,
        line width=0.8pt
    ]}
]


  \foreach \p [count=\m from 0] in {0.7,1.2,1.7, 2.5, 3.5, 4.5, 5.5, 6.5}{
    \node[point] (l1\m) at (\p, 0) {};
  }

  \node[point, color=azulejo, label={[label distance=0.7mm]270:$\scriptstyle (1,0)$}, inner sep=1.1pt] (limit) at (-.5, 0) {};

  \foreach \p [count=\m from 0] in {-0.5, 6.5}{
    \node[ipoint] (11\m) at (\p, 4.5) {};
  }

  \foreach \p [count=\m from 0] in {-0.5}{
    \node[ipoint] (21\m) at (\p, 3.9) {};
  }

  \foreach \p [count=\m from 0] in {-0.5, 5.5, 6.5}{
    \node[ipoint] (31\m) at (\p, 3.3) {};
  }

  \foreach \p [count=\m from 0] in {-0.5}{
    \node[ipoint] (41\m) at (\p, 2.7) {};
  }

  \foreach \p [count=\m from 0] in {-0.5, 4.5, 5.5, 6.5}{
    \node[ipoint] (51\m) at (\p, 2.1) {};
  }

  \foreach \p [count=\m from 0] in {-0.5}{
    \node[ipoint] (61\m) at (\p, 1.5) {};
  }

  \foreach \p [count=\m from 0] in {-0.5, 3.5 ,4.5, 5.5, 6.5}{
    \node[ipoint] (71\m) at (\p, .9) {};
  }

    \node at (-0.5,0.4) {$\vdots$};
    \node at (0.7,0.4) {$\vdots$};
    \node at (3,0.4) {$\vdots$};
    \node at (5,0.4) {$\vdots$};
    \node at (0.2, 0) {$\dots$};

     \foreach \from/\to  in {
    110/111,
    310/311,
    311/312,
    510/511,
    511/512,
    512/513,
    710/711,
    711/712,
    712/713,
    713/714}
    {
        \draw[arrow] (\from) to [bend left=0] (\to);
    }

    \foreach \from/\to  in {
    111/210,
    312/410,
    513/610}
    {
        \draw[arrow] (\from) to [bend left=5] (\to);
    }

    \foreach \from/\to  in {
    210/310,
    410/510,
    610/710}
    {
        \draw[arrow] (\from) to [bend right] (\to);
    }

    \foreach \from/\to  in {
    l10/l11,
    l11/l12,
    l12/l13,
    l13/l14,
    l14/l15,
    l15/l16,
    l16/l17}
    {
        \draw[arrow] (\from) to [bend left=20] (\to);
    }

    \draw[arrow] (l17) to [bend left=13] (limit);

    \draw[arrow, color=azulejo] (limit) to [loop] (limit);

    \node[left, font=\tiny] at (110) { $(1,1)$};
    \node[right, font=\tiny] at (111) { $f(1,1)$};
    \node[left, font=\tiny] at (210) { $f^2(1,1)$};
    \node[left, font=\tiny] at (310) { $f^4(1,1)$};

\end{tikzpicture}
\caption{Dynamical system $(Y_\alpha, f_\alpha)$ with $\alpha = (0,0, \dots)$.}
\label{fig: sistema dinamico f_alfa}
\end{figure}
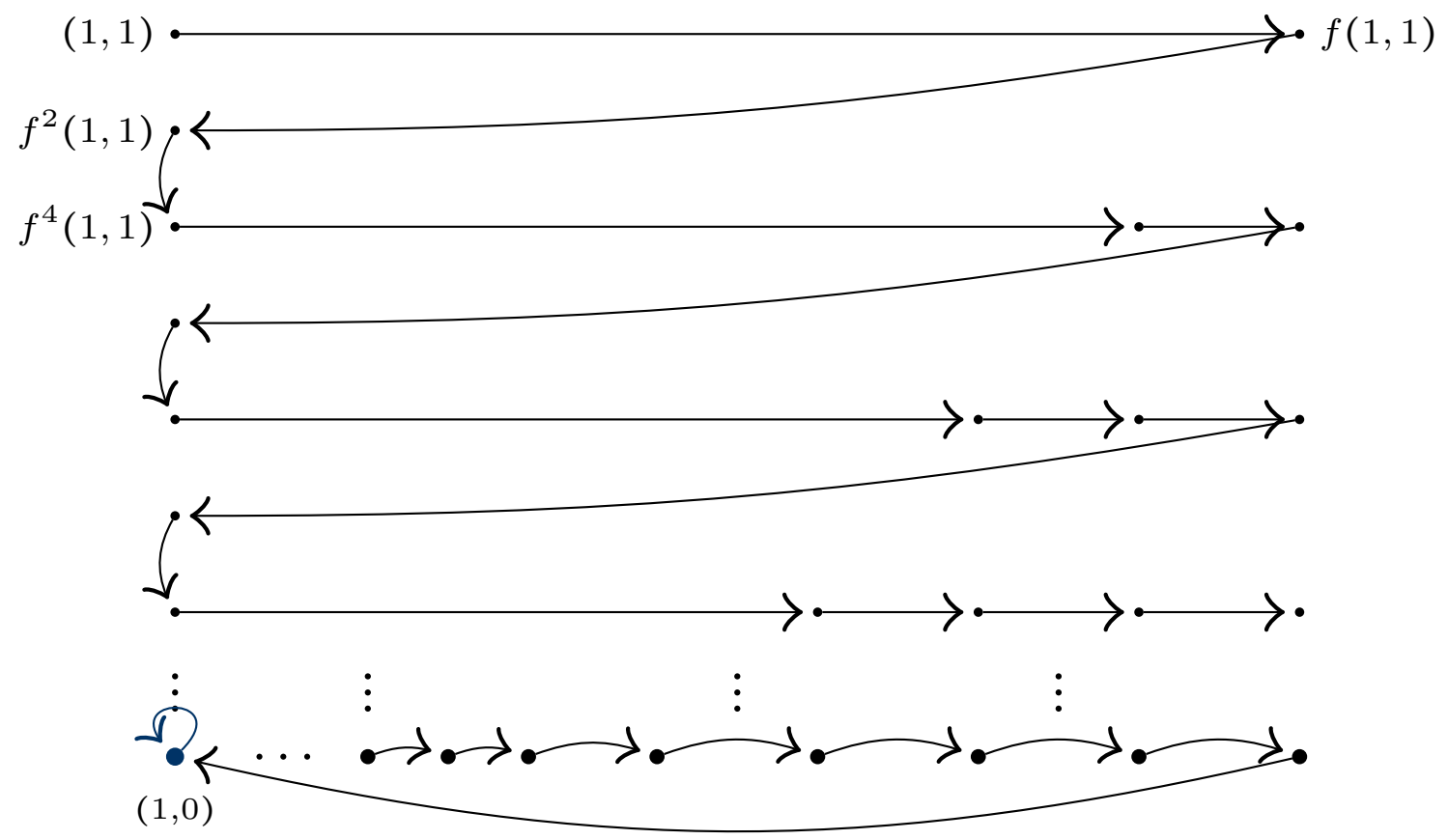
    
    We claim that $(Y_\alpha, f_\alpha)$ satisfies the desired conditions. Observe that condition \eqref{lema-item: condicion punto fijo falpha} is already satisfied.
    
    By induction on $i \in \N$ and using Lemma \ref{propiedades de G_alpha}\eqref{imagen de G_alpha}, it is easy to check that
    \[
    \orb[f_\alpha][1,1] = \bigcup_{i \in \N} M_i \times \{ \phi(i)\}. 
    \]
    Thus, the orbit of $(1,1)$ is dense by Proposition \ref{Proposicion: funcion con orbita densa}. Therefore we get \eqref{lema-item: condicion punto con orbita densa}.

    We now prove that $(Y_\alpha, f_\alpha)$ is indeed a dynamical system; that is, $f_\alpha$ is continuous. Let $(x,y) \in Y_\alpha$. We will prove continuity at point $(x,y)$. Since $(Y_\alpha)^{(1)} = W \times \{ 0 \}$ we only need to consider the case when $x \in W$ and $y = 0$. For $\epsilon > 0$, consider the ball centered at $(x,0)$ with radius $\epsilon$ given by 
    \[
    B_\epsilon(x,0) = (x - \epsilon, x + \epsilon) \times [0, \epsilon) \cap Y_\alpha.
    \]
    If $x \in (W)^{(1)} = \{ 1 \}$, since $s_\alpha$ is continuous there is $\delta_1 \in (0, \epsilon)$ such that  if $y \in (1, 1 + \delta_1) \cap W$, then $s_\alpha(y) \in (1, 1 + \epsilon)$.
    Let $\delta_2 = w_{n_0} - 1$ where $n_0$ is such that $w_n \in (1, 1 + \epsilon)$ for all $n \geq n_0$. Define $\delta = \min \{ \delta_1, \delta_2, 1/2n_0 \}$. An easy computation shows that $f_\alpha[B_\delta(1,0)] \subseteq B_\epsilon(1,0)$. 
    If $x \in I(W)$, consider $\delta \in (0, \epsilon)$ such that $(x - \delta, x + \delta) \cap W = \{ x \}$. By the definition of the sets $M_i$ and the function $f_\alpha$ it is clear that $f_\alpha[B_\delta(x, 0)] \subseteq B_\epsilon(s_\alpha(x),0)$. Thus from the above, we conclude that $f_\alpha$ is continuous.

    Now, let us observe that
    $M = W \times \{0\}$
    is a closed $f_\alpha$-invariant subset for each $\alpha \in \cantor$. By Lemma \ref{propiedades de G_alpha}\eqref{conjugacion G_alpha 2}  the restrictions of $f_\alpha$ to $M$ are  pairwise topologically nonconjugate. To see that the functions $f_\alpha$ are pairwise topologically nonconjugate, notice that any homeomorphism between $Y_\alpha$ and $Y_\beta$ must map $M$ onto itself, since $M$ is the set of limit points for any $Y_\alpha$. This implies that if $f_\alpha$ and $f_\beta$ are conjugate, then their restrictions to $M$ would also be conjugate. Thus, we conclude \eqref{lema-item: condicion M es cerrado fi} and \eqref{lema-item: condicion falpha y fbeta no conjugados}.
\end{proof}

It is worth mentioning that the dynamical systems constructed above do not exhaust all possible cases. In particular, the systems obtained here are not injective. However, in \cite[Example 4.2]{Cardinality2018} the authors provide an example of an injective dynamical system with a dense orbit.

The following is a general property on dynamical systems with a dense orbit on $\w^2 + 1$.

\begin{prop}
    \label{prop: general prop para w2 + 1}
    Let $(X, f)$ be a countable dynamical system with a dense orbit and $X$ homeomorphic to $\w^2 + 1$. Then for every $x \in X^{(1)}$ we have $|\w(x,f)| = 1$.
\end{prop}
\begin{proof}
    We will use Lemmas 3.1 and 3.2 from \cite{Cardinality2018}. Although these lemmas are stated for phase spaces of the form $\w^\alpha + 1$, they also apply to homeomorphic phase spaces. The verification is left to the reader.
    
    If $x \in X^{(2)}$ then $f(x) = x$ by \cite[Lemma 3.2(ii)]{Cardinality2018}. On the other hand, if $x \in X^{(1)} \setminus X^{(2)}$ by \cite[Lemma 3.1(v)]{Cardinality2018} we get $x$ is not periodic. Therefore $x$ is eventually periodic or the orbit of $x$ is infinite. In any case we get that the orbit of $x$ converges to the unique point of $X^{(2)}$ since $f[X^{(1)}] \subseteq X^{(1)}$ by Lemma \ref{Proposicion: funcion con orbita densa}.
\end{proof}

\section{Existence of countable dynamical systems with a dense orbit}

In this section we prove the following.

\begin{teo}
\label{OD-no-vacio}
Let $X$ be a compact metric countable space. Then $\OD(X)$ is nonempty.
\end{teo}

Notice that $\Homeo(X)$ has cardinality $2^{\aleph_0}$, hence  each  conjugacy class has $2^{\aleph_0}$ points. As $\OD(X)$ is closed under conjugacy, and the function $\varphi \in \Homeo(X)  \mapsto \varphi \circ f \circ \varphi^{-1}$ is injective for any $f \in \OD(X)$, then $\OD(X)$ has cardinality $2^{\aleph_0}$. In fact,  we will show  a stronger result (see Theorem \ref{Teorema existencia}). If $X^{(2)}$ is not empty, there are $2^{\aleph_0}$ non conjugated systems in $\OD(X)$.

The construction of an $f\in \OD(X)$ was inspired by  \cite{ChaosInCountableDS} where the author constructs homeomorphisms exhibiting distributional chaos. Although the dynamical property considered there is different from the property studied here, we adapt the method used there to obtain systems with dense orbits.

We introduce some definitions that will be used throughout this section. Given a finite subset $A \subseteq [0,1]$, we define the predecessor and successor functions  $T_A^- : A \setminus \{ \inf A \} \to A$ and $T_A^+: A \setminus \{ \sup A \} \to A$ respectively by: 
\[
T_A^-(x) = \sup(A \cap [0, x)), \ T_A^+(x) = \inf(A \cap (x, 1]).
\]
Next, we show that these definitions can be extended to the case of an infinite discrete subset $A$, under suitable conditions.

\begin{lema}
    \label{Lema: Funciones T bien definidas}
    Let $A \subseteq (0,1)$ be a discrete subset such that: 
    \begin{enumerate}[(i)]
        \item $0,1 \in A^{(1)}$.
        \item Every point $x \in A^{(1)} \setminus \{ 0,1 \}$ is a bilateral limit point of $A$.
    \end{enumerate}
    Then $\inf (A \cap (x,1]), \sup (A \cap [0,x)) \in A$ for each $x \in A$.
\end{lema}
\begin{proof}
    We will show that $\min (A \cap (x,1])$ exists. To do this, it suffices to find $y \in A$ such that $(x,y) \cap A = \emptyset$. Let
    \[
        z=\inf \{ w \in A^{(1)} : w > x \}.
    \]
    Since $z\in A^{(1)}$, it follows from (ii) that $[x,z]\cap A$ is infinite and that $z$ is the unique limit point of $[x,z]\cap A$.
    We claim that $A\cap(x,z)$ has a minimum. Indeed, if it did not, then, since it is bounded below by $x$, its infimum would be a limit point of $A$ belonging to $(x,z)$, contradicting the fact that $z$ is the unique limit point of $[x,z]\cap A$. 
    Therefore, letting $y=\min(A\cap(x,z))$ we obtain $(x,y)\cap A=\emptyset$. Hence $y=\min(A\cap(x,1])$, and thus $\min(A\cap(x,1])$ exists.

    A similar argument, considering points to the left of $x$, shows that $\max (A  \cap [0, x))$ exists.
\end{proof}

The following is straightforward and will be used later. 
\begin{lema}
    \label{Lema: funciones T}
    Let $A \subseteq [0,1]$ be a finite subset with $|A| > 1$. Then 
    $$
    T_A^-[A \setminus \{ \inf A \}] = A \setminus \{ \sup A \} \quad \text{and} \quad T_A^+[A \setminus \{ \sup A \} ] = A \setminus  \{ \inf A \}.
    $$
\end{lema}

Let $F \subseteq \R$ be a subset and $\epsilon > 0$. A set $A \subseteq F$ is said to $\epsilon$-span $F$, if for each $x \in F$ there is $y \in A$ such that $|x - y| < \epsilon$.

\begin{lema}
    \label{lema: Existencia de sistemas con orbita densa}
    Let $X$ be a compact metric countable space with $|X^{(2)}| \geq 2$. Then there exists a countable dynamical system $(Z, f)$ such that
    \begin{enumerate}[(i)]
        \item 
        \label{lema-item: Z homeomorfo a X}
        $Z$ is homeomorphic to $X$.
        \item 
        \label{lema-item: condicion sobre Z}
        $Z \subseteq [0,1]^2$, $Z^{(1)} \subseteq [0,1] \times \{ 0 \}$ and $(1,0) \in Z^{(2)}$.
        \item 
        \label{lema-item: f es inyectiva}
        $f$ is injective.
        \item 
        \label{lema-item: punto con orbita densa bajo f}
        
        The orbit of the point $(0,1)$ under $f$ is dense.

        \item 
        \label{lema-item: puntos rank 1 convergiendo y rank 2 fijos}
        If $z \in Z^{(1)}$ then  $|\w(z,f)| = 1$. Moreover, if $z \in Z^{(2)}$ then $z$ is fixed. 
        
    \end{enumerate}
\end{lema}

\begin{proof}
    By Theorem \ref{Sierpinski}, we can find a subspace $U \subseteq [0,1]$ (with the usual topology) homeomorphic to $X^{(1)}$. Moreover, we can modify $U$ appropriately so that it satisfies the following properties:
    \begin{itemize}
        \item $0, 1 \in U^{(1)}$,
        \item For each $x \in U^{(1)} \setminus \{ 0 , 1\}$, $x$ is a bilateral limit point of $U$.
    \end{itemize}
    Fix a partition $\{ A, B \}$ of $I(U)$ such that $A^{(1)} = B^{(1)} = U^{(1)}$ and every $x \in U^{(1)} \setminus \{ 0,1\}$ is a bilateral limit point of both $A$ and $B$. Furthermore, $0$ and $1$ are limit points for both $A$ and $B$. Since $A, B \subseteq [0,1]$, we can take $\{ M_{i} \}_{i \in \N}$ such that:
    \begin{itemize}
        \item $0,1 \in M_i$ and $M_i$ is finite for every $i \in \N$,
        \item $M_{2i} \setminus \{ 0,1 \} \subseteq A$, $M_{2i + 1} \setminus \{ 0,1 \} \subseteq B$ for every $i \in \N$,
        \item $M_{2i} \setminus \{0,1\}$ $\phi(2i)$-spans $A$ and $M_{2i + 1} \setminus \{ 0, 1 \}$ $\phi(2i + 1)$-spans $B$ for each $i \in \N$,
    \end{itemize}
    where $\phi: \N \to [0,1]$ is given by $\phi(i) = \frac{1}{i + 1}$.
    Consider the following
    \[
    Z = U \times \{ 0 \} \cup \bigcup_{i \in \N} M_i \times \{ \phi(i) \}.
    \]
    The space $Z$, endowed with the subspace topology, is homeomorphic to $X$. Indeed, let us observe that $Z^{(1)} = U \times \{ 0 \}$. Therefore $Z^{(1)}$ is homeomorphic to $X^{(1)}$ since $U$ is homeomorphic to $X^{(1)}$. Thus $X$ and $Z$ are homeomorphic by Corollary \ref{Corolario: homeomorfismo de puntos limites es suficiente}. Hence $Z$ satisfies \eqref{lema-item: Z homeomorfo a X} and \eqref{lema-item: condicion sobre Z}. 

    To define $f$ on $Z^{(1)}$, recall that by Lemma \ref{Lema: Funciones T bien definidas}, the functions $T^+_A$ and $T^-_B$ are well defined. Thus, for $(x,0) \in U \times \{ 0 \}$, we let
    \begin{equation}
        \label{Ec. Definicion F caso 2}
        f(x,0) = \begin{cases}
            (x,0) & \text{ if } x \in U^{(1)},\\[4pt]
            (T_{A}^+(x), 0) & \text{ if } x \in A,\\[4pt]
            (T_{B}^-(x), 0) & \text{ if } x \in B.
        \end{cases}
    \end{equation}

    Finally, to complete the definition of $f$ on the isolated points, if $(x,\phi(i)) \in M_i \times \{ \phi(i) \}$ for some $i \in \N$, we then set
    \begin{equation}
        \label{Ec. Definicion F aislados}
        f(x, \phi(i)) = \begin{cases}
            (T^+_{M_i}(x), \phi(i)) & \text{ if } i \text{ is even and } x \in M_{i}\setminus \{ 1 \},\\[4pt]
            (1, \phi(i + 1)) & \text{ if } i \text{ is even and } x = 1,\\[4pt]
            (T^-_{M_i}(x), \phi(i)) & \text{ if } i \text{ is odd and } x \in M_{i} \setminus \{  0 \},\\[4pt]
            (0, \phi(i + 1)) & \text{ if } i \text{ is odd and } x = 0.
        \end{cases}
    \end{equation}

    By the definition of $f$, it is clear that $f$ is injective. Thus, item \eqref{lema-item: f es inyectiva} is satisfied. On the other hand, since $Z^{(1)} = U \times \{  0 \}$ we get $Z^{(2)} = U^{(1)} \times \{ 0 \}$, hence by definition of $f$ we obtain that every point of $Z^{(2)}$ is fixed. Now let us observe that if $z \in Z^{(1)} \setminus Z^{(2)}$, then $z = (x, 0)$ where $x \in A \cup B$. We assume w.l.o.g. that $x \in A$. Then $\orb[f][z] = \{ ((T_A^+)^n(x), 0) \}$. Let us note that $((T_A^+)^n(x))_{n \in \N}$ is a bounded increasing sequence. Then this sequence converges to a point in $U$, because $U$ is compact. Thus the orbit of $z$ is convergent. This proves \eqref{lema-item: puntos rank 1 convergiendo y rank 2 fijos}. 
 
    We now show that $f$ is continuous at every point.
    Let $(x,0) \in U \times \{ 0 \}$ and $\epsilon > 0$. Consider the following cases:
    \begin{itemize}
        \item {\it Case 1:} $(x,0) \in I(U) \times \{ 0 \}$. Let $(y,0) = f(x,0)$. We suppose w.l.o.g. that $x \in A$. By \eqref{Ec. Definicion F caso 2}, we have that $y \in I(U)$. Let $\delta \in (0, \epsilon)$ be such that
        \[
        (x- \delta, x + \delta) \cap U = \{x\}, \ \ (y - \delta, y + \delta) \cap U = \{ y \}. 
        \]
        We claim that $f[B_\delta(x,0)] \subseteq B_\delta(y, 0)$. Indeed, if $(u, v) \in B_\delta(x,0)$, then $(u, v) = (x, 0)$ or $(u, v) \in M_i \times \{ \phi(i) \}$ for some $i$ even and $\phi(i) < \delta$. We only need to consider the last case. Since $M_i$ $\phi(i)$-spans $A$, we have that $M_i$ $\delta$-spans $A$. It follows that $x,y \in M_i$ by the choice of $\delta$. Therefore $f(u,v) = (y, \phi(i)) \in B_\epsilon(y, 0)$.

\medskip

        \item {\it Case 2:} $(x,0) \in U^{(1)} \times \{ 0 \}$. From \eqref{Ec. Definicion F caso 2} we have $f(x,0) = (x,0)$. First, we assume that $x \notin \{0, 1\}$. 
        By definition of $A$ and $B$, we can take $y_0 \in A$, $z_0 \in B$ such that $x- \epsilon < z_0 < x < y_0 < x + \epsilon$. Now, let $\delta \in (0, \epsilon)$ be such that
        \begin{align*}
            & \ \ \ x - \epsilon < z_0 < x - \delta < x < x + \delta < y_0 < x+\epsilon,\\
            & \delta < \min \Big\{ \frac{1}{2}|x - y_0|, \frac{1}{2}|x - z_0|, \epsilon - |x- y_0|, \epsilon - |x - z_0| \Big\}.
        \end{align*}
        We claim that $f[B_\delta(x,0)] \subseteq B_\epsilon(x,0)$. Indeed, let $(u, v) \in B_\delta(x,0)$.
        If $v = 0$, then $u \in U^{(1)}$ or $u \in I(U)$. Note that $f(u, 0) = (u, 0) \in B_\epsilon(x,0)$ whenever $u \in U^{(1)}$. Thus, we consider the remaining case. We assume w.l.o.g. that $u \in A$, thus we get
        \[
        u < T_{A}^+(u) < x, \text{ whenever } u < x \quad \text{ or } \quad u < T_{A}^+(u) \leq y_0, \text{ whenever } u > x.
        \]
        Therefore, we conclude that $f(u, 0) = (T_A^+(u), 0) \in B_\epsilon(x,0)$.
    
        On the other hand, if $v \neq 0$, then $(u,v) \in M_i \times \{ \phi(i) \}$ for some $i \in \N$. We suppose w.l.o.g that $i$ is even. Then we have $u \in M_i \setminus \{ 0,1 \} \subseteq A$. If $u = 1$, it is easy to see that $f(1, \phi(i)) \in B_\epsilon(x,0)$. Therefore we consider $u \neq 1$. Now let us note that $u < y_0 - \delta$ and $y_0 + \delta < x + \epsilon$. Hence $T_{M_i}^+(u) \leq y_0 + \delta$ since $M_i$ $\delta$-spans $A$ (because $\phi(i) < \delta$). Thus $f(u,v) = (T_{M_i}^+(u), v) \in B_\epsilon(x,0)$.

        \medskip
        If $x = 0$ or $1$, a similar argument applies. 
        
        \end{itemize}

        By induction on $i$ and Lemma \ref{Lema: funciones T}, we have 
        \[
            \orb[f][0, 1] = \bigcup_{i \in \N} M_i \times \{ \phi(i) \}.
        \]
        Hence, $(Z, f)$ is a dynamical system with a dense orbit.
\end{proof}

\bigskip

\noindent {\em Proof of Theorem \ref{OD-no-vacio}.}
We proceed by cases, according to the Cantor--Bendixson rank of $X$. If $X^{(2)} = \emptyset$, then by Theorem \ref{Sierpinski}, $X$ is homeomorphic to $k \cdot \w + 1$ for some $k \in \N$, and the result follows from the dynamical systems constructed in Section \ref{Section: case kw + 1}. Otherwise, if $X^{(2)} \neq \emptyset$, the result follows from Lemmas \ref{lema: Existencia sistemas dinamicos en w2 + 1} and \ref{lema: Existencia de sistemas con orbita densa}.
\qed

\subsection{ \texorpdfstring{$\OD(X)$}{OD(X)} is nowhere dense}
In this  short section we show that $\OD(X)$ is topologically small, namely, nowhere dense with respect to the compact-open topology  on $C(X,X)$.

\begin{lema}
    \label{lema: nunca denso}
    Let $X$ be a compact zero\--dimensional metric  space, $A \subsetneq X$ nonempty, and let
    \[
        \mathcal{I}(A) = \{ f \in \C(X, X) : f(A) \subseteq A \}.
    \]
    If $A$ is closed, then $\mathcal{I}(A)$ is closed. If furthermore $X \setminus A$ is dense, then $\mathcal{I}(A)$ is nowhere dense.
\end{lema}
\begin{proof}
    Let $A \subseteq X$ be a closed set. We show that $\mathcal{I}(A)$ is closed. Indeed, let $e_x: \C(X,X) \to X$ be the evaluation function (that is $e_x(f) = f(x)$). Let us note that 
    \[
        \mathcal{I}(A) = \bigcap_{x \in A} e_{x}^{-1}[A].
    \]
    Since $A$ is closed and $e_x$ is continuous for each $x$, it follows that $\mathcal{I}(A)$ is closed.

    Now assume that $X \setminus A$ is dense. We show that $\mathcal{I}(A)$  has empty interior. Let $f \in \mathcal{I}(A)$ and $U$ be an open set with $f \in U$. We can assume that $U$ is a basic open set in the compact\--open topology, that is, there exists $n \in \N$ and there are $K_i$ and $U_i$ compact and open sets (respectively) of $X$, for $i = 1, \dots, n$, such that 
    \[
    U = \{ g \in \C(X, X) : g[K_i] \subseteq  U_i \text{ for each $i = 1, \dots, n$} \}.
    \]
    We have the following cases:

    \begin{enumerate}[(i)]
        \item $ A \cap (\bigcup_{i = 1}^n K_i)$ \textit{is empty}. Let $V$ be a clopen set  such that $V \cap K_i = \emptyset$ for each $i$, and $A \cap V \neq \emptyset$. Take $y \in X \setminus A$ and consider $h: X \to X$ given by
        \[
            h(x) = \begin{cases}
                y & \text{ if } x \in V,\\
                f(x) & \text{ if } x \notin V.
            \end{cases}
        \]
        Then $h$ is continuous and $h \in U$ but $h \notin \mathcal{I}(A)$.

        \item $A \cap (\bigcup_{i = 1}^n K_i)$ \textit{is nonempty}. Let $z \in A\cap (\bigcup_{i = 1}^n K_i)$, and let $M = \{ i \in \{1,\dots, n\} : z \in K_i \}$. Since $X \setminus A$ is dense, there is $y \in X \setminus A$  with $y \in \bigcap_{i \in M} U_i$ (let us note that the intersection is nonempty because $f(z)$ belongs to this intersection).Choose a clopen set $V \subseteq X$ with $z \in V$ such that $V \cap K_i = \emptyset$ for $i \notin M$ (if $M=\{1,2,\ldots, n\}$ take any clopen set containing $z$). Define $h : X \to X$ by
        \[
        h(x) = \begin{cases}
                y & \text{ if } x \in V,\\
                f(x) & \text{ if } x \notin V.
            \end{cases}
        \]
        Let us note that $h$ is continuous, and $h \notin \mathcal{I}(A)$ because $h(z) = y \notin A$. We claim that $h \in U$. Indeed, fix $i = 1, \dots, n$ and $x \in K_i$. We show that $h(x) \subseteq U_i$. Consider the following cases: \textit{(i)} $i \in M$. If $x \notin V$, then $h(x) = f(x) \in U_i$. Otherwise, if $x \in V$, then $h(x) = y \in U_i$ by the choice of $y$. \textit{(ii)} $i \notin M$. By the choice of $V$, we get $V \cap K_i = \emptyset$. In particular, $x \notin V$, hence $h(x) = f(x) \in U_i$.
    \end{enumerate}

    From the above, we obtain $f \in U \not\subseteq \mathcal{I}(A)$. 
    Since $f$ and $U$ were arbitrary, we conclude that $\mathcal{I}(A)$ has empty interior.  
\end{proof}

\begin{teo}
    Let $X$ be a compact metric countable space. Then $\OD(X)$ is nowhere dense.
\end{teo}
\begin{proof}
    Indeed, let us note that $\OD(X) \subseteq \mathcal{I}(X^{(1)})$ by Proposition \ref{Proposicion: funcion con orbita densa}. Hence the conclusion follows immediately from Lemma \ref{lema: nunca denso}.
\end{proof}
\section{Gluing countable dynamical systems with a dense orbit}
\label{Section: pegado}

In this section we introduce a gluing procedure for countable dynamical systems with a dense orbit. To the best of our knowledge, this gluing construction is new. It remains unclear whether its compatibility conditions can be weakened.

We first describe how to glue the orbit of one point to the orbit of another under suitable compatibility conditions. This construction is then extended to gluing two  dynamical systems. The results obtained in this section will play a crucial role in the proof of the main theorem of this article.

\subsection{Gluing orbits}

Let $X_1, X_2$ be two infinite sets. For $i = 1,2$, let $f_i : X_i \to X_i$ and $x_i \in X_i$  such that the orbit of $x_i$ under $f_i$ is infinite.  Let $Y = \orb[f_1][x_1] \cup \orb[f_2][x_2]$.
We give conditions under which there exists a map $g:Y\to Y$ such that $\orb[g][x_1] = Y$ and $g$  is \textit{consistent} with both $f_1$ and $f_2$. We call $g$ the {\it gluing function of the orbit of $x_2$ to $x_1$}.

We say that the orbit of $x_2$ {\it is compatible with the orbit of} $x_1$ if the following conditions hold: 
\begin{itemize}
    \item[(C1)] $x_2 \in \orb[f_1][x_1]$,
    \item[(C2)] $K = \orb[f_1][x_1] \cap \orb[f_2][x_2]$ is infinite. Let $(n_i)_{i \in \N}$ and $(m_i)_{i \in \N}$ be the increasing enumerations of
    \[
        \{ n_i : i \in \N \} = \{ n \in \N : f^n_1(x_1) \in K \} \ \text{ and } \ \{ m_i:  i \in \N \} =  \{ n \in \N : f^n_2(x_2) \in K \},
    \]
    respectively.
    \item[(C3)] $f_1^{n_i}(x_1) = f_2^{m_i}(x_2)$ for all $i \in \N$,
    \label{compatibilidad - condicion sobre n,m}
    \item[(C4)] $n_{2i + 1} = n_{2i} + 1$ and $m_{2i + 2} = m_{2i + 1} + 1$ for all $i \in \N$.
\end{itemize}

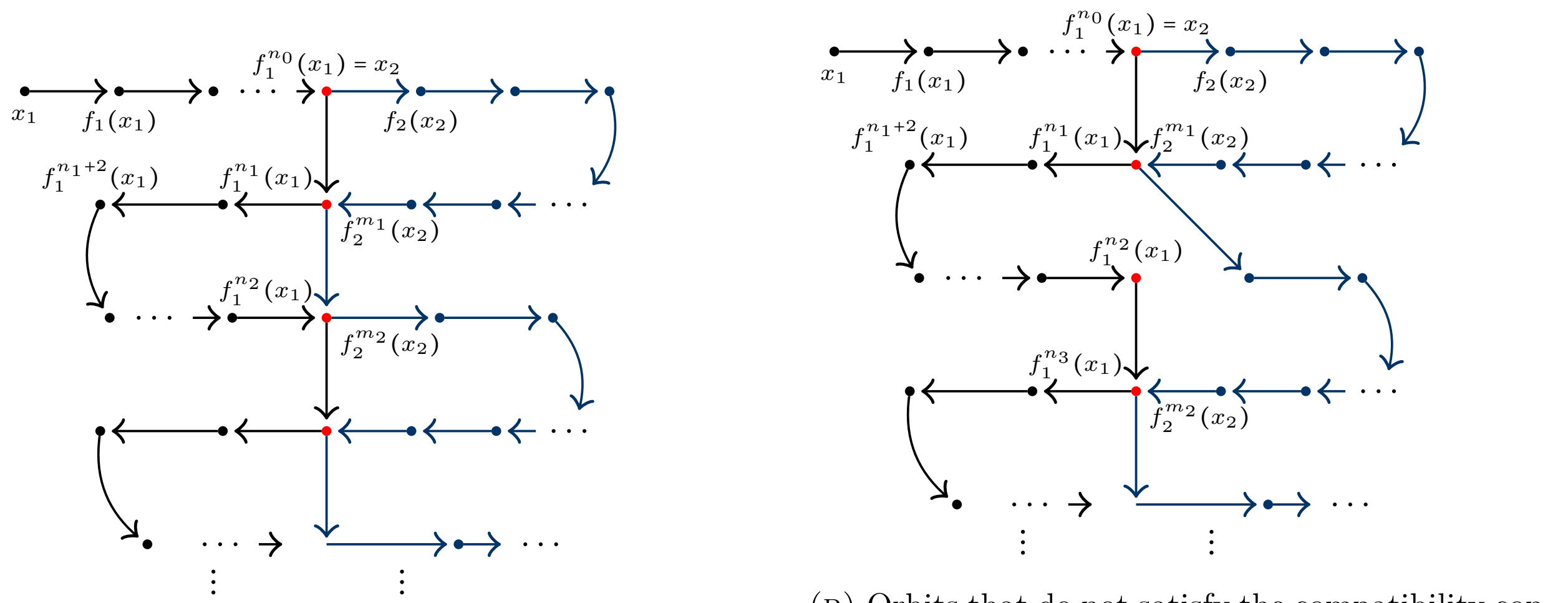
\begin{figure}[H]
    \centering

    \begin{subfigure}{0.48\textwidth}
        \begin{tikzpicture}[
    every label/.append style={font=\footnotesize},
    scale=1
  ]

    \tikzstyle{point}=[draw,shape=circle, inner sep=.9pt, fill]
    
    \tikzstyle{ipoint}=[draw,shape=circle, inner sep=.5pt, fill]

    \tikzstyle{arrow} = [
    line width=0.7pt,
    -{Classical TikZ Rightarrow[
        scale=1.4,
        sep=2pt,
        length=1mm,
        line width=0.9pt
    ]}
    ]

    \definecolor{azulejo}{RGB}{0,51,102}

  \foreach \p [count=\m from 0] in {1,2,3}{
    \node[point] (\m) at (\p, 0) {};
  }

  \node (3) at (3.5,0) {$\dots$};
  
  \node          [point, color=red] (4) at (4.2, 0) {};
  \node[point, color=red] (5) at (4.2, -1.2) {};

  \foreach \p [count=\m from 6] in {3.1, 1.8}{
    \node[point] (\m) at (\p, -1.2) {};
  }

  \node[point] (8) at (1.9,-2.4) {};
  \node (9) at (2.4,-2.4) {$\dots$};

  \node[point] (10) at (3.2,-2.4) {};

  \node[point, color=red] (11) at (4.2, -2.4) {};
  \node[point, color=red] (12) at (4.2, -3.6) {};
  
  \foreach \p [count=\m from 13] in {3.1, 1.8}{
    \node[point] (\m) at (\p, -3.6) {};
  }

  \node[point] (15) at (2.3,-4.8) {};
  \node (16) at (3.1,-4.8) {$\dots$};

    \foreach \from/\to in {
        0/1,
        1/2,
        3/4,
        4/5,
        5/6,
        6/7,
        9/10,
        10/11,
        11/12,
        12/13,
        13/14}
    {
        \draw[arrow] (\from) to (\to);
    }

    \foreach \from/\to in {
    7/8,
    14/15}
    {
        \draw[arrow] (\from) to [bend right] (\to);
    }

    \draw[arrow] (16) to (3.8,-4.8);


    \foreach \p [count=\m from 17] in {5.2, 6.2, 7.2}
    {
        \node[point, color=azulejo] (\m) at (\p, 0) {};
    }

    \node (20) at (6.8,-1.2) {$\dots$};
    \node[point, color=azulejo] (21) at (6,-1.2) {};
    \node[point, color=azulejo] (22) at (5.1,-1.2) {};

    \foreach \p [count=\m from 23] in {5.4, 6.6}
    {
        \node[point, color=azulejo] (\m) at (\p, -2.4) {};
    }

    \node (25) at (6.8,-3.6) {$\dots$};
    \node[point, color=azulejo] (26) at (6,-3.6) {};
    \node[point, color=azulejo] (27) at (5.1,-3.6) {};

    \node[point, color=azulejo] (28) at (5.6, -4.8) {};
    \node (29) at (6.5, -4.8) {$\dots$};

    \node at (3,-5.2) {$\vdots$};

    \node at (5,-5.2) {$\vdots$};

    \foreach \from/\to in {
        4/17,
        17/18,
        18/19,
        20/21,
        21/22,
        22/5,
        5/11,
        11/23,
        23/24,
        25/26,
        26/27,
        27/12,
        28/29}
    {
        \draw[arrow, color=azulejo] (\from) to (\to);
    }

    \foreach \from/\to in {
        19/20,
        24/25}{
        \draw[arrow, color=azulejo] (\from) to [bend left] (\to);
    }

    \draw[arrow, color=azulejo] (12) to (4.2,-4.8);
    \draw[arrow, color=azulejo] (4.2, -4.8) to (28);

    \node[below=1.7pt, font=\tiny] at (0) {$x_1$};
    \node[below=1.7pt, font=\tiny] at (1) {$f_1(x_1)$};
    \node[above, font=\tiny] at (4) {$f_1^{n_0}(x_1) = x_2$};

    \node[below=1.7pt, font=\tiny] at (17) {$f_2(x_2)$};

    \node[above left, font=\tiny] at (5) {$f_1^{n_1}(x_1)$};

    \node[below right, font=\tiny] at (5) {$f_2^{m_1}(x_2)$};

    \node[above, font=\tiny] at (7) {$f_1^{n_1 + 2}(x_1)$};

    \node[above left, font=\tiny] at (11) {$f_1^{n_2}(x_1)$};

    \node[below right, font=\tiny] at (11) {$f_2^{m_2}(x_2)$};
\end{tikzpicture}
        \caption{Orbits under suitable conditions}
        \label{figure: pegado orbitas 1}
    \end{subfigure}
    \hfill
    \begin{subfigure}{0.48\textwidth}
        \begin{tikzpicture}[
    every label/.append style={font=\footnotesize},
    scale=1
  ]

    \tikzstyle{point}=[draw,shape=circle, inner sep=.9pt, fill]
    
    \tikzstyle{ipoint}=[draw,shape=circle, inner sep=.5pt, fill]

    \tikzstyle{arrow} = [
    line width=0.7pt,
    -{Classical TikZ Rightarrow[
        scale=1.4,
        sep=2pt,
        length=1mm,
        line width=0.9pt
    ]}
    ]

    \definecolor{azulejo}{RGB}{0,51,102}

  \foreach \p [count=\m from 0] in {1,2,3}{
    \node[point] (\m) at (\p, 0) {};
  }

  \node (3) at (3.5,0) {$\dots$};
  
  \node          [point, color=red] (4) at (4.2, 0) {};
  \node[point, color=red] (5) at (4.2, -1.2) {};

  \foreach \p [count=\m from 6] in {3.1, 1.8}{
    \node[point] (\m) at (\p, -1.2) {};
  }

  \node[point] (8) at (1.9,-2.4) {};
  \node (9) at (2.4,-2.4) {$\dots$};

  \node[point] (10) at (3.2,-2.4) {};

  \node[point, color=red] (11) at (4.2, -2.4) {};
  \node[point, color=red] (12) at (4.2, -3.6) {};
  
  \foreach \p [count=\m from 13] in {3.1, 1.8}{
    \node[point] (\m) at (\p, -3.6) {};
  }

  \node[point] (15) at (2.3,-4.8) {};
  \node (16) at (3.1,-4.8) {$\dots$};

    \foreach \from/\to in {
        0/1,
        1/2,
        3/4,
        4/5,
        5/6,
        6/7,
        9/10,
        10/11,
        11/12,
        12/13,
        13/14}
    {
        \draw[arrow] (\from) to (\to);
    }

    \foreach \from/\to in {
    7/8,
    14/15}
    {
        \draw[arrow] (\from) to [bend right] (\to);
    }

    \draw[arrow] (16) to (3.8,-4.8);


    \foreach \p [count=\m from 17] in {5.2, 6.2, 7.2}
    {
        \node[point, color=azulejo] (\m) at (\p, 0) {};
    }

    \node (20) at (6.8,-1.2) {$\dots$};
    \node[point, color=azulejo] (21) at (6,-1.2) {};
    \node[point, color=azulejo] (22) at (5.1,-1.2) {};

    \foreach \p [count=\m from 23] in {5.4, 6.6}
    {
        \node[point, color=azulejo] (\m) at (\p, -2.4) {};
    }

    \node (25) at (6.8,-3.6) {$\dots$};
    \node[point, color=azulejo] (26) at (6,-3.6) {};
    \node[point, color=azulejo] (27) at (5.1,-3.6) {};

    \node[point, color=azulejo] (28) at (5.6, -4.8) {};
    \node (29) at (6.5, -4.8) {$\dots$};

    \node at (3,-5.2) {$\vdots$};

    \node at (5,-5.2) {$\vdots$};

    \foreach \from/\to in {
        4/17,
        17/18,
        18/19,
        20/21,
        21/22,
        22/5,
        5/23,
        23/24,
        25/26,
        26/27,
        27/12,
        28/29}
    {
        \draw[arrow, color=azulejo] (\from) to (\to);
    }

    \foreach \from/\to in {
        19/20,
        24/25}{
        \draw[arrow, color=azulejo] (\from) to [bend left] (\to);
    }

    \draw[arrow, color=azulejo] (12) to (4.2,-4.8);
    \draw[arrow, color=azulejo] (4.2, -4.8) to (28);

    \node[below=1.7pt, font=\tiny] at (0) {$x_1$};
    \node[below=1.7pt, font=\tiny] at (1) {$f_1(x_1)$};
    \node[above, font=\tiny] at (4) {$f_1^{n_0}(x_1) = x_2$};

    \node[below=1.7pt, font=\tiny] at (17) {$f_2(x_2)$};

    \node[above left, font=\tiny] at (5) {$f_1^{n_1}(x_1)$};

    \node[above right, font=\tiny] at (5) {$f_2^{m_1}(x_2)$};

    \node[above, font=\tiny] at (7) {$f_1^{n_1 + 2}(x_1)$};

    \node[above, font=\tiny] at (11) {$f_1^{n_2}(x_1)$};

    \node[above left, font=\tiny] at (12) {$f_1^{n_3}(x_1)$};
    
    \node[below right, font=\tiny] at (12) {$f_2^{m_2}(x_2)$};
\end{tikzpicture}
        \caption{Orbits that do not satisfy the compatibility conditions}
        \label{figure: pegado orbitas 2}
    \end{subfigure}

    \caption{Compatibility conditions for gluing orbits}
    \label{figure: pegado de orbitas}
\end{figure}

The compatibility conditions are designed to ensure that the two orbits can be concatenated without skipping orbit segments. Figure \ref{figure: pegado de orbitas} illustrates the idea. In Figure \ref{figure: pegado orbitas 1}, the successive intersection points (red points) are arranged so that one can alternate between the two orbits while preserving the intervening orbit segments. At each intersection, the construction switches from one orbit to the other in a controlled manner.

Starting from $x_1$, we follow its orbit until it reaches the first intersection point, $f_1^{n_0}(x_1)=x_2$. At this point, we continue along the orbit of $x_2$ until reaching the next intersection point, $f_2^{m_1}(x_2)=f_1^{n_1}(x_1)$. We then proceed along the orbit of $x_1$ again and repeat this process alternately. Thus, the resulting orbit alternately follows finite segments of the two original orbits, with the intersection points determining where the switches occur.

By contrast, Figure  \ref{figure: pegado orbitas 2} illustrates a configuration that is not suitable for our construction.
In this case, switching between the orbits at an intersection point (at the point $(f_1^{n_1}(x_1) = f_2^{m_1}(x_2))$) would require jumping over a segment of one of the original orbits, as indicated by the diagonal arrow. 
Consequently, part of the orbit would no longer be represented in the glued orbit. Conditions (C3) and (C4) prevent this behavior.

We now define the gluing function $g: \orb[f_1][x_1] \cup \orb[f_2][x_2] \to \orb[f_1][x_1] \cup \orb[f_2][x_2]$ as follows.
\begin{equation}
\label{definicion funcion de pegado 1}
g(z) = \begin{cases}
    f_1(z) & \text{ if } z \in \orb[f_1][x_1] \setminus \orb[f_2][x_2],\\[3pt]
    f_2(z) & \text{ if } z \in \orb[f_2][x_2] \setminus \orb[f_1][x_1],\\[3pt]
    f_2(z) & \text{ if } z \in \orb[f_1][x_1] \cap \orb[f_2][x_2] \text{ and } z = f_1^{n_i}(x_1) \text{ for some even } i,\\[3pt]
    f_1(z) & \text{ if } z \in \orb[f_1][x_1] \cap \orb[f_2][x_2] \text{ and } z = f_1^{n_i}(x_1) \text{ for some odd } i.
\end{cases}
\end{equation}

We now verify that $g$ glues the two orbits as intended.

\begin{lema}
    \label{Lema funcion de pegado}
     Let $X_1, X_2$ be two infinite sets, $f_i : X_i \to X_i$ and $x_i \in X_i$ with infinite orbit under $f_i$ for $i = 1,2$. Assume that the orbit of $x_2$ is compatible with the orbit of $x_1$, and let $g$ be the gluing function. Then $\orb[g][x_1] = \orb[f_1][x_1] \cup \orb[f_2][x_2]$.
\end{lema}
\begin{proof}
    By (C3) and the definition of $g$ it is clear that $\orb[g][x_1] \subseteq \orb[f_1][x_1] \cup \orb[f_2][x_2]$. Now we proceed to show that $\orb[f_1][x_1] \subseteq \orb[g][x_1]$ and $\orb[f_2][x_2] \subseteq \orb[g][x_1]$.

    We introduce a sequence $(s_i)_{i\in\N}$ specifying the iterates of $g$ that lie in the intersection $\orb[f_1][x_1] \cap \orb[f_2][x_2]$. More precisely, we describe the sequence $(s_i)_{i\in\N}$ so that
    
    \begin{equation}
    \label{Ecuacion: sucesion m_i}
    g^{s_i}(x_1) = f_1^{n_i}(x_1) = f_2^{m_i}(x_2) \quad \text{ for all } i \in \N.
    \end{equation}
    
    We set $s_0 = n_0$, and $s_1 = n_0 + m_1$. In general, for $i > 1$ 
    \[
    \begin{cases}
        s_i = n_i + m_{i - 1} - (i - 1) \text{ whenever } i \text{ is even,}\\
        s_i = n_{i -1} + m_i - (i - 1) \text{ whenever }  i \text{ is odd}.
    \end{cases}
    \]
    We will verify, by induction on $i$, that the sequence $(s_i)_{i \in \N}$ actually satisfies equation \eqref{Ecuacion: sucesion m_i}. First, let us note that 
    \begin{equation}
        \label{Eq: resta de s_i}
        \begin{cases}
            s_{i + 1} - s_i = n_{i + 1} - n_i & \text{ whenever } i \text{ is odd},\\
            s_{i + 1} - s_i = m_{i + 1} - m_i & \text{ whenever } i \text{ is even}.
        \end{cases}
    \end{equation}
    
    Second, observe that $g^j(x_1) = f_1^j(x_1)$ for all $j \leq n_0$ due to minimality of $n_0$.
    Third, assume that \eqref{Ecuacion: sucesion m_i} is satisfied for $i$. If  $i$ is even, let $t = m_{i + 1} - m_i$. By definition of $g$ and \eqref{Ecuacion: sucesion m_i} we have 
    \[
    g(g^{s_i}(x_1)) = f_2(g^{s_i}(x_1)) = f_2(f_2^{m_i}(x_2)).
    \]
    Thus, for each $0 < r < t$ we obtain $f_2^{m_i + r}(x_2) \notin \orb[f_1][x_1]$, which implies that 
    \begin{equation}
    \label{Ecuacion: propiedad g y s_i}
        g^{s_i + r}(x_1) = f_2^{m_i + r}(x_2) \text{ for } 0 \leq r < t = m_{i + 1} - m_i.
    \end{equation}
    Hence, for $r = t - 1$ we have
    \[
        g^{s_i + t}(x_1) = g(g^{s_i + t - 1}(x_1)) = g(f_2^{m_i + t - 1}(x_2)) = f_2(f_2^{m_i + t - 1}(x_2)) = f_2^{m_i + t}(x_2)
    \]
    Now let us observe that $m_i + t = m_{i + 1}$ and $s_i + t = s_{i + 1}$ by \eqref{Eq: resta de s_i}. Therefore $g^{s_{i + 1}}(x_1) = f_2^{m_{i + 1}}(x_2) = f_1^{n_{i + 1}}(x_1)$.

    If $i$ is odd and taking $t = n_{i + 1} - n_i$ (and using a similar argument to the case in which $i$ even) we obtain 
    \begin{equation}
        \label{Ecuacion propiedad g y s_i -2}
        g^{s_i + r}(x_1) = f_1^{n_i + r}(x_1) \text{ for } 0 \leq r \leq t.
    \end{equation}    
    Thus, for $r = t$ and using \eqref{Eq: resta de s_i} we have $g^{s_{i + 1}}(x_1) = f_1^{n_{i + 1}}(x_1) = f_2^{m_{i + 1}}(x_2)$. From the above cases, we conclude the proof of \eqref{Ecuacion: sucesion m_i}.
    
    We next prove that $\orb[f_2][x_2] \subseteq \orb[g][x_1]$. Indeed, let $j \in \N$. Then $j \in [m_i, m_{i + 1}]$ for some $i \in \N$. If $i$ is even, by equation \eqref{Ecuacion: propiedad g y s_i} we obtain that $f_2^j(x_2) \in \orb[g][x_1]$.
    On the other hand, if $i$ is odd, let us note that $[m_i,m_{i + 1}] = \{ m_i, m_{i + 1} \}$ by condition (C4). Therefore, by equation \eqref{Ecuacion: sucesion m_i} the conclusion follows.

    A completely similar argument, using \eqref{Ecuacion propiedad g y s_i -2} instead of \eqref{Ecuacion: propiedad g y s_i} and that $g^j(x_1) = f_1^j(x_1)$ for all $j \leq n_0$, we conclude $\orb[f_1][x_1] \subseteq \orb[g][x_1]$.
\end{proof}

\subsection{Gluing dynamical systems}
{$\;$}
\label{section: pegado de sistemas}

In this subsection we extend the orbit gluing construction to  the dynamical systems themselves. Under suitable compatibility assumptions, we show that two countable dynamical systems with dense orbits can be glued into a new system that preserves both continuity and the existence of a dense orbit.

Throughout this section, all compact metric countable spaces are regarded as subspaces of $[0,2] \times [0,1]$. Hence, for compact metric countable spaces $X_1$ and $X_2$, the union $X_1 \cup X_2$ carries the subspace topology inherited from $[0,2] \times [0,1]$.

\begin{lema}
\label{pegando sistemas dinamicos}
Let $(X_1, f_1)$ and $(X_2, f_2)$ be two countable dynamical systems with dense orbits and $x_i \in X_i$ for $i = 1,2$. Let $X = X_1 \cup X_2$ and suppose that the following statements hold:
\begin{enumerate}[(i)]
    \item The orbit of $x_i$ under $f_i$ is dense in $X_i$ for $i = 1,2$.
    \item The orbit of $x_2$ is compatible with the orbit of $x_1$.
    \item
    \label{Prop: Condicion igualdad puntos aislados}
    $I(X) = I(X_1) \cup I(X_2)$.
    \item 
    \label{Prop: Condicion igualdad de funciones}
    $f_1(x) = f_2(x)$ for all $x \in X_1^{(1)} \cap X_2^{(1)}$.
\end{enumerate}
Let  $g$ be the gluing function of the orbit of $x_2$ to $x_1$ and $f_1 \oplus f_2 : X \to X$  be given by
\begin{equation}
\label{Ecuacion: funcion pegado de sistemas}
    (f_1 \oplus f_2) (x) = \begin{cases}
        g(x) & \text{ if } x \in I(X),\\
        f_1(x) & \text{ if } x \in X_1^{(1)},\\
        f_2(x) & \text{ if } x \in X_2^{(1)}.
    \end{cases}
\end{equation}
Then $(X, f_1 \oplus f_2)$ is a countable dynamical system with a dense orbit.
\end{lema}
\begin{proof}
    Notice that $f_1 \oplus f_2$ is well defined by \eqref{Prop: Condicion igualdad puntos aislados} and \eqref{Prop: Condicion igualdad de funciones}. 
    From Proposition \ref{Proposicion: funcion con orbita densa}, we have that $I(X_i) = \orb[f_i][x_i]$ for $i = 1,2$. On the other hand, by \eqref{Prop: Condicion igualdad puntos aislados} and Lemma \ref{Lema funcion de pegado}, we conclude that the orbit of $x_1$ under $f_1 \oplus f_2$ is $I(X)$. It remains to prove that $f_1 \oplus f_2$ is a continuous function. 
    
    To do this, take $z \in X^{(1)}$ and  a sequence $(z_n)_{n \in \N}$ in $X$ such that $(z_n)_{n \in \N}$ converges to $z$. First let us assume that $(z_n)_{n \in \N} \subseteq I(X)$. For $i = 1,2$, define
    \[
    A_i = \{ n \in \N: g(z_n) = f_i(z_n) \text{ and } z_n \in X_i \}.
    \]
    By definition of $g$ (see \eqref{definicion funcion de pegado 1}) we have $A_1 \cup A_2 = \N$. If both $A_1$ and $A_2$ are infinite, then $z \in X_1 \cap X_2$, and therefore $z \in X_1^{(1)} \cap X_2^{(1)}$ by \eqref{Prop: Condicion igualdad puntos aislados}. In this case, the subsequences $(g(z_n))_{n \in A_i}$ converge to $f_i(z)$ for $i = 1,2$, and condition \eqref{Prop: Condicion igualdad de funciones} implies that $(g(z_n))_{n \in \N}$ converges to $f_1(z) = f_2(z)$. 
    If instead one of the sets $A_{i}$ is infinite and the other is finite, then the continuity of the corresponding function $f_i$ shows that $(g(z_n))_{n \in \N}$ converges to $f_{i}(z)$, and hence $((f_1 \oplus f_2)(z_n))_{n \in \N}$ converges to $(f_1 \oplus f_2)(z)$.

    Second, we assume that $(z_n)_{n \in \N} \subseteq X^{(1)}$. We argue as in the previous case. For $i = 1,2$, set 
    \[
    B_i = \{ n \in \N : z_n \in X_i \}.
    \]
    If both $B_1$ and $B_2$ are infinite, then $z\in X_1^{(1)} \cap X_2^{(1)}$. Since each $f_i$ is continuous, the subsequences $(f_i(z_n))_{n \in B_i}$ converge to $f_i(z)$, and by definition of $f_1 \oplus f_2$ together with \eqref{Prop: Condicion igualdad de funciones}, it follows that the sequence $((f_1 \oplus f_2)(z_n))_{n \in \N}$ converges to $(f_1 \oplus f_2)(z)$. 
    If only one of the sets $B_i$ is infinite, then the continuity of the corresponding function $f_i$ ensures the convergence of $(f_i(z_n))_{n \in \N}$, and therefore $((f_1 \oplus f_2)(z_n))_{n \in \N}$ converges to $(f_1 \oplus f_2)(z)$.
\end{proof}

Whenever two countable dynamical systems $(X_1, f_1)$ and $(X_2, f_2)$ satisfy the hypotheses of Lemma \ref{pegando sistemas dinamicos}, we say that \textit{$(X_1, f_1)$ can be glued to $(X_2, f_2)$}. 

The following lemma will be used to construct nonconjugate dynamical systems. 

\begin{lema}
    \label{distinguiendo sistemas dinamicos con el pegado}
    Let $(X_1,f_1), (X_2, f_2)$ and $(Y, h)$ be three countable dynamical systems with a dense orbit such that:
    \begin{enumerate}[(i)]
        \item $(Y, h)$ can be glued to both $(X_1, f_1)$ and $(X_2, f_2)$.
        \item
        \label{Lema: h inyectiva}
        $h$ is one to one.
        \item There exists an infinite closed subset $M \subseteq X_1 \cap X_2$ such that $M$ is strongly invariant under $f_1$ and $f_2$; that is, $f_i[M] = M$ for $i = 1,2$.
        \item 
        \label{Lema: union y conjugacion - condicion punto fijo}
        Both $(X_1,f_1)$ and $(X_2,f_2)$  have a unique fixed point $z$. Furthermore, $z \in M \cap Y$ and satisfies
        \[
            \bigcup_{n \in \N} f_{i}^{-n}[z] = M, \ \text{ for } i = 1,2.
        \]

        \item
        $(M, f_1)$ and $(M, f_2)$ are nonconjugate.
    \end{enumerate}
    Then $(Y \cup X_1, h \oplus f_1)$ and $(Y \cup X_2, h \oplus f_2)$ are also nonconjugate.
\end{lema}
\begin{proof}
    For $i = 1,2$, define $g_i = h \oplus f_i$. Suppose that $(Y \cup X_1, g_1)$ and $(Y \cup X_2, g_2)$ are conjugate, and let $\varphi \in \Homeo(Y \cup X_1, Y \cup X_2)$ be such that $\varphi \circ g_1 = g_2 \circ \varphi$. Note that $M \subseteq X_i^{(1)}$ because each point of $M$ is eventually periodic and no isolated point can be eventually periodic.
    Now we will show that $(M, f_1)$ and $(M, f_2)$ are also conjugate. More precisely, we will prove that $\varphi[M] = M$; hence, the restriction $\varphi\restriction_M$ is the required conjugacy between $(M, f_1)$ and $(M, f_2)$ because $g_i = f_i$ on $M$.
    
    Let us observe that $g_i(z) = f_i(z) = z$ since $z \in X_i$ for $i \in \{ 1,2\}$.

    \begin{claimInProof}
        \label{Claim: punto fijo bajo phi}
        $\varphi(z)=z$. Suppose for contradiction that $\varphi(z) \neq z$. Since $\varphi(z) = \varphi(g_1(z)) = g_2(\varphi(z))$, it follows that $\varphi(z)$ is a fixed point under $g_2$.
        Now we show that  $\varphi(z) \notin X_2$. Indeed, if $\varphi(z) \in X_2$, then by definition of \(g_2\) we would have $g_2(\varphi(z))=f_2(\varphi(z))$. Since $\varphi(z)$ is fixed under $g_2$, it follows that $f_2(\varphi(z)) = \varphi(z)$, which contradicts \eqref{Lema: union y conjugacion - condicion punto fijo}. Therefore $\varphi(z) \in Y \setminus X_2$. 

        On the other hand, by the definition of $g_2$ and \eqref{Ecuacion: funcion pegado de sistemas} we have $h^{-1}[\varphi(z)] \subseteq g^{-1}_2[\varphi(z)]$ (notice that $z$ is a limit point, therefore, $\varphi(z)$ is also a limit point). Moreover, since $\varphi(z) \notin X_2$ we also have $g^{-1}_2[\varphi(z)] \subseteq h^{-1}[\varphi(z)]$. Thus $g^{-1}_2[\varphi(z)] = h^{-1}[\varphi(z)]$.

        Since \(\varphi\) is a conjugacy, we have $\varphi[g_1^{-1}[z]] = g_2^{-1}[\varphi(z)]$. Moreover, by the definition of \(g_1\) we obtain $f_1^{-1}[z] \subseteq g_1^{-1}[z]$. Since $|f_1^{-1}[z]| > 1$ it follows that $|g_1^{-1}[z]| > 1$. Hence, 
        $$
            |\varphi[g_1^{-1}[z]]| = |h^{-1}[\varphi(z)]| > 1.
        $$
        Which is a contradiction since $h$ is injective. Therefore, $\varphi(z) = z$.
        \qed \ Claim \ref{Claim: punto fijo bajo phi}.
    \end{claimInProof}

    \begin{claimInProof}
        \label{Claim: contenencia preimagen g_i}
        \textit{For all $n \in \N$ and $i \in \{ 1,2 \}$,
        \[
        g_i^{-n}[z] \subseteq X_i.
        \]}
         
        First note that every point of \(g_i^{-n}[z]\) is a limit point. In fact, recall that $g_i(z) = z$. This implies that those points are eventually mapped to the fixed point \(z\), and hence are eventually periodic. Now notice that $h = f_i$ on $Y^{(1)} \cap X_i^{(1)}$ because $h$ and $f_i$ can be glued (see the hypothesis \eqref{Prop: Condicion igualdad de funciones} of Lemma \ref{pegando sistemas dinamicos}). In particular, $h(z) = f_i(z) = z$.
        
        We proceed by induction on \(n\).
        For the base case \(n = 1\), let \(y \in g_i^{-1}[z]\) and suppose that \(y \notin X_i\). By the definition of \(g_i\) (see \eqref{Ecuacion: funcion pegado de sistemas}), it follows that \(g_i(y) = h(y)\), and hence \(h(y) = z\). Since \(y \notin X_i\), we have \(y \neq z\), which implies that \(h\) is not injective, contradicting \eqref{Lema: h inyectiva}.

        Assume that \(g_i^{-n}[z] \subseteq X_i\). Let \(y \in g_i^{-(n+1)}[z]\) and suppose that \(y \notin X_i\). Then \(g_i(y) \in g_i^{-n}[z]\), so by the inductive hypothesis \(g_i(y) \in X_i\). On the other hand, since \(y \notin X_i\), we have \(g_i(y) = h(y)\), and thus \(h(y) \in X_i \cap Y\). Moreover,
        \[
            g_i^{n+1}(y) = h^{n+1}(y) = z,
        \]
        and again \(y \neq z\). This shows that \(h\) is not injective, a contradiction.
        \qed \ Claim \ref{Claim: contenencia preimagen g_i}.
    \end{claimInProof}
    
    \begin{claimInProof}
        \label{Claim: preimagenes iguales bajo g y f}
        \textit{
        For all \(n \in \N\) and \(i \in \{1,2\}\),}
        \[
            g_i^{-n}[z] = f_i^{-n}[z].
        \]
        From definition of $g_i$ (see \eqref{Ecuacion: funcion pegado de sistemas}) it is clear that $f_i^{-n}[z] \subseteq g_i^{-n}[z]$. For the reverse inclusion let \(y \in g_i^{-n}[z]\). By Claim \ref{Claim: contenencia preimagen g_i}, we have $y \in X_i$. Hence, by definition of $g_i$, it follows that $g_i(y) = f_i(y)$, and therefore $g_i^{n}(y) = f_i^{n}(y)$. Since $g_i^n(y) = z$, we conclude that $f_i^n(y) = z$. This proves the reverse inclusion.
        \qed \ Claim \ref{Claim: preimagenes iguales bajo g y f}
    \end{claimInProof}

    \medskip
    
    We now have everything needed to show that $\varphi[M] = M$. By part \eqref{Lema: union y conjugacion - condicion punto fijo} and Claim \ref{Claim: preimagenes iguales bajo g y f}, we obtain
    \[
    \bigcup_{n \in \N} g_i^{-n}[z] = M, \qquad \text{for } i = 1,2.
    \]
    On the other hand, by Lemma \ref{Lema: preimagen conjugacion} and Claim \ref{Claim: punto fijo bajo phi}, we have
    \[
    \varphi[g_1^{-n}[z]] = g_2^{-n}[\varphi(z)] = g_2^{-n}[z], \qquad \text{for all } n \in \N.
    \]
    Taking the union over $n \in \N$ on both sides, and using the identity above, we conclude that $\varphi[M] = M$.
\end{proof}

\section{The main theorem}

In this section we prove the main theorem of the paper. The result provides a dichotomy for dynamical systems with dense orbit, depending on the second Cantor–Bendixson derivative of the phase space $X$.

We prove  that if the second derivative of $X$ is empty, then all dynamical systems with a dense orbit on $X$ are topologically conjugate; that is, all such systems are essentially the ones constructed in Section \ref{Section: case kw + 1}. On the other hand, when the second derivative of $X$ is nonempty, we show that there exist continuum many pairwise nonconjugate such systems. The proof relies on the constructions developed in the previous sections, together with the gluing technique introduced earlier.

\begin{teo}
    \label{Teorema existencia}
    Let $X$ be a countable compact metric space. Then the following statements hold.
    \begin{enumerate}[(i)]
        \item If $X^{(2)} = \emptyset$, then any two functions in $\OD(X)$ are topologically conjugate.
        
        \item \label{Teorema existencia - 2} 
        If $X^{(2)} \neq \emptyset$, then there are $2^{\aleph_0}$ topologically nonconjugate functions in $\OD(X)$. 
    \end{enumerate}
\end{teo}
\begin{proof}
        {\it (i)} Let $f_1, f_2 \in \OD(X)$ and $x_1, x_2$ be the points whose orbits  are dense under $f_1$ and $f_2$, respectively. As $X^{(2)} = \emptyset$, it follows that $X^{(1)}$ is finite. Set $X^{(1)} = \{ z_1, z_2, \dots, z_k \}$. Now, let us note that $\w(x_1, f_1) = \w(x_2, f_2) = X^{(1)}$, since the orbits of $x_1$ and $x_2$ are dense. 
        By Lemma \ref{caracterizacion convergencia en conjuntos omega limites}, for each $j = 1, \dots, k$, there exists a unique $0 \leq r_j < k$ such that by letting $A_j = k \N + r_j$, the sequence \[
        (f_1^n(x_1))_{n \in A_j} \ \ \text{converges to } z_j.
        \] 
        Furthermore (using again Lemma \ref{caracterizacion convergencia en conjuntos omega limites}) we get that, for all $j = 1, \dots, k$,
    \[
        (f_1^n(x_1))_{n \in A_j} \text{ converges } \iff (f_2^n(x_2))_{n \in A_j} \text{ converges}.
    \]
    Thus, the following function is well defined: Define $\varphi: X \to X$ by
    \[
    \varphi(x) = \begin{cases}
        f_2^n(x_2) & \text{ if } x = f_1^n(x_1) \text{ for some } n \in \N,\\
        \lim_{n \in A_{j}} f_2^n(x_2) & \text{ if } x = z_j \text{ for some } 1 \leq j \leq k.
    \end{cases}
    \]
    By Lemma \ref{caracterizacion convergencia en conjuntos omega limites} we conclude  that $\varphi$ is a continuous function. 
    Using Proposition \ref{Proposicion: funcion con orbita densa} it is easy to check that $\varphi$ is bijective. 
    Finally, by continuity and the definition of $\varphi$, we obtain $\varphi \circ f_1 = f_2 \circ \varphi$. This concludes the proof.

    {\it (ii)} If $|X^{(2)}| = 1$, then $X$ is homeomorphic to $\w^2 + 1$. Therefore the  conclusion follows from Lemma  \ref{lema: Existencia sistemas dinamicos en w2 + 1} and Lemma \ref{OrbitaDensaBajoHomeomorfismos}. Now assume that $|X^{(2)}| > 1$. Using Lemma \ref{OrbitaDensaBajoHomeomorfismos}, it is sufficient to construct, for each $\alpha \in \cantor$, a space $Z_\alpha$ homeomorphic to $X$ and a function $f_\alpha \in \OD(Z_\alpha)$ such that $f_\alpha$ and $f_\beta$ are nonconjugate whenever $\alpha \neq \beta$. 
    
    Let $\{ (Y_\alpha, f_\alpha) \}_{\alpha \in \cantor}$ and $(Z,f)$ be the dynamical systems obtained from Lemmas \ref{lema: Existencia sistemas dinamicos en w2 + 1} and \ref{lema: Existencia de sistemas con orbita densa} respectively. 

    \begin{claimInProof}
        \label{claim: homeo X}
        \textit{$Z_\alpha := Z \cup Y_\alpha$ is homeomorphic to $Z$}. Indeed, by Lemma \ref{lema: Existencia sistemas dinamicos en w2 + 1}\eqref{lema-item: condicion sobre Yalpha} and Lemma \ref{lema: Existencia de sistemas con orbita densa}\eqref{lema-item: condicion sobre Z}, it follows that the second Cantor-Bendixson derivative of $Y_\alpha$ is contained in the second derivative of $Z$.
        Thus by Proposition \ref{union homeomorfa a uniendo} we conclude the claim.
        \qed \ Claim \ref{claim: homeo X}
    \end{claimInProof}

    \begin{claimInProof}
        \label{claim: pegado de sistemas}
        \textit{$(Z, f)$ can be glued to $(Y_\alpha, f_\alpha)$ for each $\alpha \in \cantor$.} Fix $\alpha \in \cantor$. Let $x_1 = (0,1)$ and $x_2 = (1,1)$. From Lemma \ref{lema: Existencia sistemas dinamicos en w2 + 1}\eqref{lema-item: condicion punto con orbita densa} and \ref{lema: Existencia de sistemas con orbita densa}\eqref{lema-item: punto con orbita densa bajo f} we have that $x_1$ and $x_2$ are the points whose orbit under $f$ and $f_\alpha$, respectively, are dense. We now verify the hypotheses of Lemma \ref{pegando sistemas dinamicos}. 
        
        First, we are going to check that the orbit of $x_2$ is compatible with the orbit of $x_1$. Since $x_2 \in \orb[f][x_1]$, condition (C1) holds. Now let us note $I(Z) \cap I(Y_\alpha) = \{ (1, \varphi(i)) : i \in \N \}$ and the orbits of $x_1$ and $x_2$ traverse $I(Z) \cap I(Y_\alpha)$ following the order $(1, \varphi(i))$, thus we get (C3). On the other hand, by definition of $f$ and $f_\alpha$ we obtain 
        \[
            f(1, \phi(2i)) = (1, \phi(2i + 1)) \quad \text{ and } \quad
            f_\alpha(1, \phi(2i + 1)) = (1, \phi(2i + 2)),
        \]
        therefore we have (C4).

        Second, it is clear that $I(Z_\alpha) \subseteq I(Z) \cup I(Y_\alpha)$ by definition of $Z_\alpha$. On the other hand, since Lemmas \ref{lema: Existencia sistemas dinamicos en w2 + 1}\eqref{lema-item: condicion sobre Yalpha} and \ref{lema: Existencia de sistemas con orbita densa}\eqref{lema-item: condicion sobre Z} every point of $Z$ (and $Y_\alpha$) is isolated if and only if its second coordinate is not zero. The above implies that $I(Z) \cup I(Y_\alpha) \subseteq I(Z_\alpha)$. Third, $(Z)^{(1)} \cap (Y_\alpha)^{(1)} = \{ (1,0) \}$ and $f_\alpha(1,0) = f(1,0) = (1,0)$. Thus, the claim follows from Lemma \ref{pegando sistemas dinamicos}.
        \qed \ Claim \ref{claim: pegado de sistemas}
    \end{claimInProof}

    Recall $Z_\alpha = Z \cup Y_\alpha$ and let $j_\alpha = f \oplus f_\alpha$. By Claim \ref{claim: pegado de sistemas} $(Z_\alpha, j_\alpha)$ is a dynamical system for each $\alpha \in \cantor$. With the following claim we finish the proof.

    \begin{claimInProof}
        \label{claim: no conjugados - final}
        \textit{For each $\alpha \neq \beta$, $(Z_\alpha, j_\alpha)$ and $(Z_\beta, j_\beta)$ are not conjugate}. Indeed, we will verify that the hypotheses of Lemma \ref{distinguiendo sistemas dinamicos con el pegado} are satisfied by taking  $ (X_1, f_1) = (Y_\alpha, f_\alpha)$, $(X_2, f_2 ) = (Y_\beta, f_\beta)$ and  $(Y, h) = (Z, f)$. We shall check each hypothesis of Lemma \ref{distinguiendo sistemas dinamicos con el pegado}:
        \begin{enumerate}[(i)]
            \item follows from Claim \ref{claim: pegado de sistemas},
            
            \item follows from Lemma \ref{lema: Existencia de sistemas con orbita densa} \eqref{lema-item: f es inyectiva},
            
            \item follows from Lemma \ref{lema: Existencia sistemas dinamicos en w2 + 1}\eqref{lema-item: condicion M es cerrado fi},
    
            \item follows from Lemma \ref{lema: Existencia sistemas dinamicos en w2 + 1}\eqref{lema-item: condicion punto fijo falpha}, and
    
            \item follows from Lemma \ref{lema: Existencia sistemas dinamicos en w2 + 1}\eqref{lema-item: condicion falpha y fbeta no conjugados}.
        \end{enumerate}

        Consequently, all the hypotheses of Lemma \ref{distinguiendo sistemas dinamicos con el pegado} are fulfilled. Therefore, $(Z_\alpha, j_\alpha)$ and $(Z_\beta, j_\beta)$ are nonconjugate.
        \qed \ Claim \ref{claim: no conjugados - final}
    \end{claimInProof}

    Finally, each $Z_\alpha$ is homeomorphic to $X$ because $Z_\alpha$ is homeomorphic to $Z$ (by Claim \ref{claim: homeo X}), and $Z$ is homeomorphic to $X$ by Lemma \ref{lema: Existencia de sistemas con orbita densa}\eqref{lema-item: Z homeomorfo a X}. 
  \end{proof}

\section{On the Ellis semigroup}

In this final short section we present a result about the Ellis semigroup $E(X,f)$ of a countable dynamical system $(X,f)$  having a dense orbit.  In \cite{Cardinality2018,GRC2021,GarciaRodriguezTomita} investigate when $E(X,f)$ is homeomorphic to $X$.  We show that the examples we constructed are of that type. On the other hand,  \cite[Theorem 3.11]{iteratesOfDN} shows that whenever every point in $X'$ is periodic, either all functions in $E(X,f)$ are continuous or only the $f^n$ are continuous for $n\in \N$. Under the assumption that $|\w(x,f)| = 1$ for each $x \in X^{(1)}$, we characterize when each of these two alternatives occurs.

We present some brief preliminaries on the Ellis semigroup needed to show our result.
For a countable dynamical system $(X,f)$ the \emph{Ellis semigroup} of $(X,f)$ is
\[
    E(X,f)=\overline{\{f^n:n\in\mathbb{N}\}}\subseteq X^X,
\]
where the closure is taken in the topology of pointwise convergence and $f^0 = \mathrm{Id}_X$. It is a compact semigroup under composition of functions.

An \emph{ultrafilter} on $\mathbb{N}$ is a family
$p\subseteq\mathcal{P}(\mathbb{N})$ such that:
\textit{(i)} $\mathbb{N}\in p$, \textit{(ii)} $\varnothing\notin p$, \textit{(iii)} $p$ is closed under finite intersections and supersets, and \textit{(iv)} for every $A\subseteq\mathbb{N}$,
either $A\in p$ or $\mathbb{N}\setminus A\in p$.
An ultrafilter is \emph{free} if it contains no finite set. We denote the set of all free ultrafilters on $\mathbb{N}$ by $\mathbb{N}^\ast$.

For any ultrafilter $p$ and a sequence $(x_n)_{n \in \N}$ in $X$, we write $x=p\--\lim_{n \to \infty}x_n$ if for every neighborhood $U$ of $x$ we have that $\{n\in\mathbb{N}:x_n\in U\}\in p$. Such a point exists and is unique because $X$ is compact and Hausdorff.
The \emph{$p$-iterate} of $f$ is the function $f^p :  X\to X$ defined by
\[
    f^p(x) = p \--\lim_{n \to \infty}f^n(x).
\]
The Ellis semigroup admits the representation
\[
    E(X,f)=\{f^p:p \text{ is an ultrafilter on } \N \}.
\]
We denote with $E^\ast(X,f)$ the set of nontrivial elements of Ellis semigroup, that is, 

\[
E^\ast(X,f) = E(X,f) \setminus \{ f^n : n \in \N \}.
\]

\bigskip 

The following results are well\--known and easy  to prove. 

\begin{lema}
    \label{lema: p y n iteradas conmutan}
    Let $(X,f)$ be a dynamical system, $n \in \N$ and $p \in \N^\ast$. Then $f^p \circ f^n = f^n \circ f^p$.
\end{lema}

\begin{lema}
    \label{lema: p iterada es limite}
    Let $(X,f)$ be a countable dynamical system, $x \in X$ and $p \in \N^\ast$. If $\orb[f][x] \subseteq I(X)$ is infinite then $f^p(x) \in X^{(1)}$.
\end{lema}

For a countable dynamical system $(X,f)$ let
\[
    C(f) = \bigg\{ y \in X : \text{ there is } x \in X^{(1)} \text{ such that } \lim_{n \to \infty} f^n(x) = y \bigg\}.
\]

\begin{teo}
\label{Ellis}
Let $(X,f)$ be a countable dynamical system with a dense orbit such that $|\w(x,f)| = 1$ for each $x \in X^{(1)}$. Then the following statements holds:
\begin{enumerate}[(i)]
    \item
    \label{teo: ellis - cf 2}
    If $|C(f)| \geq 2$ then  $E^\ast(X,f)$ contains only discontinuous functions.
    
    \item 
    \label{teo: ellis - cf 1}
    If $|C(f)| = 1$ then $E^\ast(X,f)$ contains only continuous functions.
\end{enumerate}
\end{teo}

\begin{proof}
Let $w \in X$ be the point with dense orbit. Since the orbit of $w$ is infinite, it is easy to check that $f^p(w) \neq f^n(w)$ for any $n \in \N$ and $p \in \N^\ast$. In particular $f^n \neq f^p$, which implies that $E^\ast(X,f) = \{ f^p : p \in \N^\ast \}$.
\begin{claimInProof}
    \label{claim: igualdad de p iteradas}
    \textit{If $p \in \N^\ast$ and $z = f^p(w)$ then}
    \[
        f^p(x) = \begin{cases}
            f^n(z) & \text{ if } x = f^n(w) \text{ for some } n,\\
            \lim_{n \to \infty} f^n(x) & \text{ if } x \in X^{(1)}.
        \end{cases}
    \]
For  $x \in X^{(1)}$, the assumption $|\w(x,f)| = 1$ and the definition of $p$\--iterate imply that $f^p(x) = \lim_{n \to \infty} f^n(x)$, thus we only need to prove that $f^p(x) = f^n(z)$ if $x = f^n(w)$. Indeed let $x = f^n(w)$, then using Lemma \ref{lema: p y n iteradas conmutan} we have
    \[
    f^p(x) = f^p(f^n(w)) = f^n(f^p(w)) = f^n(z).
    \]
    This proves the claim.
    \qed \ Claim \ref{claim: igualdad de p iteradas}
\end{claimInProof}

We prove  \textit{(i)}, that is, we show that $f^p$ is discontinuous for every $p \in \N^\ast$. Let $z_1 = f^p(w)$ and $y_1 = \lim_{n \to \infty} f^n(z_1)$. By Lemma \ref{lema: p iterada es limite} we have $z_1 \in X^{(1)}$ which implies that $y_1 \in C(f)$. Since $|C(f)| \geq 2$ there is $y_2 \in C(f)$ with $y_1 \neq y_2$. As $y_2 \in C(f)$ there is $z_2 \in X^{(1)}$ such that $y_2 = \lim_{n \to \infty} f^n(z_2)$. We will show that $f^p$ is discontinuous at $z_2$.

Since $I(X)$ is dense in $X$, there is a sequence $(x_i)_{i \in \N} \subseteq I(X)$ converging  to $z_2$. Let $m_i \in \N$ be such that $x_i = f^{m_i}(w)$ for each $i$. Since $x_2 \notin I(X)$ the sequence $(m_i)_{i \in \N}$ converges to $\infty$. Hence $\lim_{i \to \infty} f^p(x_i) = \lim_{i \to \infty} f^{m_i}(z_1) = y_1$ by  Claim \ref{claim: igualdad de p iteradas}. But $f^p(z_2) = \lim_{n \to \infty} f^n(z_2) = y_2$ and $y_2 \neq y_1$, thus $f^p$ is discontinuous at $z_2$.

To prove \textit{(ii)} let us observe that by Claim \ref{claim: igualdad de p iteradas}, the equality $C(f) = \{ y \}$ implies $f^p(x) = y$ for any $x \in X^{(1)}$. Let $z \in X^{(1)}$, we check that $f^p$ is continuous at $z$. Since $f^p$ is constant on $X^{(1)}$, it suffices to prove that for every sequence $(x_i)_{i \in \N} \subseteq I(X)$ that converges to $z$ we have $f^p(x_i)$ converges to $f^p(z) = y$. Let $m_i \in \N$ with $f^{m_i}(w) = x_i$. Then $f^p(x_i) = f^{m_i}(f^p(w))$. But $f^p(w) \in X^{(1)}$ by Lemma \ref{lema: p iterada es limite}, thus $\lim_{i \to \infty} f^p(x_i) = \lim_{i \to \infty} f^{m_i}(f^p(w)) = y$. Hence $f^p$ is continuous at $z$.
\end{proof}

Finally, let $H: E(X,f)\to X$ be given by $H(f^p)=f^p(w)$. Clearly $H$ is continuous. It is easy to see that if $w$ has a dense orbit, then $H$ is onto. Under the hypotheses of the previous theorem it follows from Claim \ref{claim: igualdad de p iteradas} that $H$ is also injective and thus $E(X,f)$ is homeomorphic to $X$, a result that also follows from \cite[Theorem 4.2]{GRC2021}.

\begin{ex}
Let $X$ be a compact metric countable space with $|X^{(2)}| \geq 2$ and $\{ (Z_\alpha, j_\alpha) : \alpha \in \cantor \}$ be the family of dynamical systems with dense orbit constructed in Theorem \ref{Teorema existencia} for $X$. Then the following statements holds: 
    \begin{enumerate}[(i)]
        \item $E(Z_\alpha, j_\alpha)$ is homeomorphic to $X$.
        \item $E^\ast(Z_\alpha, j_\alpha)$ contains only discontinuous functions.
    \end{enumerate}
\end{ex}
\begin{proof}
    \textit{(i)} Recall that the dynamical systems $(Z_\alpha, j_\alpha)$ are obtained by gluing two dynamical systems. More precisely, each is obtained by gluing of the dynamical system $(Z,f)$ and $(Y_\alpha, f_\alpha)$ obtained from Lemma \ref{lema: Existencia de sistemas con orbita densa} and \ref{lema: Existencia sistemas dinamicos en w2 + 1} respectively. Now recall that these dynamical systems satisfy $|\w(z,f)| = 1$ at every limit point $z$ (for $(Z,f)$ see Lemma \ref{lema: Existencia de sistemas con orbita densa}\eqref{lema-item: puntos rank 1 convergiendo y rank 2 fijos} and for $(Y_\alpha, f_\alpha)$ see Proposition \ref{prop: general prop para w2 + 1}). Hence the dynamical systems $(Z_\alpha,j_\alpha)$ have the same property (see \eqref{Ecuacion: funcion pegado de sistemas}). Therefore by \cite[Theorem 4.2]{GRC2021} we obtain $E(Z_\alpha, j_\alpha)$ is homeomorphic to $Z_\alpha$. By Claim \ref{claim: homeo X} we have $Z_\alpha$ is homeomorphic to $Z$. Furthermore from Lemma \ref{lema: Existencia de sistemas con orbita densa}\eqref{lema-item: Z homeomorfo a X} $Z$ is homeomorphic to $X$, thus we conclude \textit{(i)}.

    To prove \textit{(ii)} let us observe that every point $x \in Z^{(2)}$ is fixed under $f$ by Lemma \ref{lema: Existencia de sistemas con orbita densa}\eqref{lema-item: puntos rank 1 convergiendo y rank 2 fijos}. Hence each point $x \in Z^{(2)}$ is fixed under $j_{\alpha}$ since $j_\alpha = f$ on $Z^{(2)}$. By the above we obtain $Z^{(2)} \subseteq C(j_\alpha)$. On the other hand, recall $X$ is homeomorphic to $Z$ by Lemma \ref{lema: Existencia de sistemas con orbita densa}\eqref{lema-item: Z homeomorfo a X}, therefore $|X^{(2)}| = |Z^{(2)}|$. Hence $2 \leq |Z^{(2)}| \leq |C(j_\alpha)|$. Finally the conclusion follows from Theorem \ref{Ellis}\eqref{teo: ellis - cf 2}.
\end{proof}

\section*{Declaration on the use of AI}
The mathematical results, proofs, and arguments presented in this article were developed by the authors without the aid of artificial intelligence tools. Generative AI was used solely for English grammar checking, manuscript polishing, and correcting typographical errors in the final version.

\bigskip

\noindent{\bf Acknowledgment.}
We thank  Universidad Industrial de Santander for  partial support.

\end{document}